\documentclass[preprint,12pt]{elsarticle}
\usepackage{graphicx}
\usepackage{amsmath,amssymb,amsthm,stmaryrd,natbib}
\usepackage[all]{xy}
\xyoption{arc}

\numberwithin{equation}{section}

\newtheorem{thm}[equation]{Theorem}
\newtheorem{cor}[equation]{Corollary}
\newtheorem{prop}[equation]{Proposition}
\newtheorem{lem}[equation]{Lemma}
\theoremstyle{definition}
\newtheorem{defn}[equation]{Definition}
\theoremstyle{remark}
\newtheorem{rmk}[equation]{Remark}
\newtheorem{exam}[equation]{Example}

\def\co{\colon\thinspace}
\newcommand{\mb}[1]{\mathbb{#1}}
\newcommand{\mf}[1]{\mathfrak{#1}}

\newcommand{\overto}{\mathop\rightarrow}

\newcommand{\into}{\mathop\hookrightarrow}

\newcommand{\TAF}{{\rm TAF}}
\newcommand{\taf}{{\rm taf}}

\newcommand{\End}{{\rm End}}

\newcommand{\Spec}{{\rm Spec}}

\newcommand{\tmf}{{\rm tmf}}

\newcommand{\TMF}{{\rm TMF}}
\newcommand{\Lie}{{\rm Lie}}

\newcommand{\Tr}{{\rm Tr}}
\newcommand{\Nm}{{\rm N}}
\newcommand{\PGL}{{\rm PGL}}
\newcommand{\Cl}{{\rm Cl}}
\newcommand{\mass}[1]{\left|{#1}\right|}
\newcommand{\card}[1]{\#\left\{#1\right\}}
\newcommand{\lsym}[2]{\left\{\frac{#1}{#2}\right\}}

\newcommand{\Zos}{\mb Z[1/6]}

\newcommand{\eilm}[1]{\ensuremath{{\mb H} #1}}

\newcommand{\comp}[1]{\ensuremath{#1^\wedge}}

\newcommand{\pow}[1]{\left\llbracket{#1}\right\rrbracket}
\newcommand{\BPP}[1]{{\rm BP}\langle{#1}\rangle}

\newcommand{\dfour}{\sigma_4}
\newcommand{\dsix}{\sigma_6}
\newcommand{\U}{W}
\newcommand{\Z}{\mathbb Z}
\newcommand{\F}{\mathbb F}
\newcommand{\sqrtdelta}{\sqrt{\Delta}}
\newcommand{\EO}{\mathit{EO}}

\newcommand{\toda}[1]{\langle#1\rangle}

\newcommand{\xym}[1]{
\vskip 0.7pc
\centerline{\xymatrix{#1}}
\vskip 0.7pc
}

\newcommand{\xynm}[1]
{\vskip 0.7pc
\centerline{\xy #1 \endxy}
\vskip 0.7pc}

\journal{Advances in Mathematics}

\begin{document}

\begin{frontmatter}

\title{Automorphic forms and cohomology theories on Shimura curves of
  small discriminant}

\author[uva]{Michael Hill}
\ead{mikehill@virginia.edu}
\author[umn]{Tyler Lawson\fnref{grant}}
\ead{tlawson@math.umn.edu}
\fntext[grant]{Partially supported by NSF grant 0805833.}

\address[uva]{Department of Mathematics,
University of Virginia\\
Charlottesville, VA 22904-4137}
\address[umn]{Department of Mathematics,
University of Minnesota\\
Minneapolis, MN 55455}

\begin{abstract}
  We apply Lurie's theorem to produce spectra associated to
  $1$-dimensional formal group laws on the Shimura curves of
  discriminants $6$, $10$, and $14$.  We compute rings of
  automorphic forms on these curves and the homotopy of the associated
  spectra.  At $p=3$, we find that the curve of discriminant $10$
  recovers much the same as the topological modular forms spectrum,
  and the curve of discriminant $14$ gives rise to a model of a
  truncated Brown-Peterson spectrum as an $E_\infty$ ring spectrum.
\end{abstract}


\end{frontmatter}

\section{Introduction}

The theory $\TMF$ of topological modular forms associates a cohomology
theory to the canonical $1$-dimensional formal group law on the
(compactified) moduli of elliptic curves.  The chromatic data of this
spectrum is determined by the height stratification on the moduli of
elliptic curves, and the homotopy groups can be computed by means of a
small Hopf algebroid associated to Weierstrass curves \cite{bauer}.

With newly available machinery this can be generalized to similar
moduli.  Associated to a discriminant $N$ which is a product of an
even number of distinct primes, there is a Shimura curve that
parametrizes abelian surfaces with action of a division algebra.
These are similar in spirit to the moduli of abelian varieties studied
previously in topology \cite{taf,absurf}.  Shimura curves carry
$1$-dimensional $p$-divisible groups which are split summands of the
$p$-divisible group of the associated abelian surface, and Lurie's
theorem allows these to give rise to associated $E_\infty$ ring
spectra analogous to $\TMF$.  See Section~\ref{sec:spectra} for
details.

In addition, there are further quotients of these Shimura curves by
certain Atkin-Lehner involutions $w_d$ for $d | N$.  Given a choice of
square root of $\pm d$ in $\mb Z_p$, we can lift this involution to an
involution on the associated $p$-divisible group, and produce
spectra associated to quotients of the Shimura curve.

Our goal in this paper is a computational exploration of spectra
associated to these curves for the three smallest valid discriminants:
$6$, $10$, and $14$.  These computations rely on a computation of
rings of automorphic forms on these curves over $\mb Z[1/N]$ in the
absence of anything approaching a Weierstrass equation.  The
associated cohomology theories, like $\TMF$, only detect information
in chromatic layers $0$ through $2$.  One might ultimately hope that a
library of techniques at low chromatic levels will allow inroads to
chromatic levels $3$ and beyond.

At discriminant $6$, we find that the ring of automorphic forms away
from $6$ is the ring
\[
\mb Z[1/6,U,V,W]/(U^4 + 3V^6 + 3W^2).
\]
See Section~\ref{sec:disc6}. This recovers a result of
Baba-Granath in characteristic zero \cite{babagranath}. The homotopy of
the associated spectrum consists of this ring together with a Serre dual
portion that is concentrated in negative degrees and degree $3$.  We
also exhibit a connective version, with the above ring as its ring of
homotopy groups.

At discriminant $14$, we obtain an application to Brown-Peterson
spectra.  Baker defined generalizations of Brown-Peterson and
Johnson-Wilson spectra with homotopy groups that are quotients of
${\rm BP}_*$ by regular sequence of generators other than the Araki or
Hazewinkel generators \cite{baker}, and Strickland made use of similar
constructions in studying products on $MU$-modules \cite{strickland}.
We show that there exists such a generalized truncated Brown-Peterson
spectrum $\BPP{2}'$ admits an $E_\infty$ structure at the prime $3$.
Associated to a prime $p$, these ring spectra ${\rm BP}$ and $\BPP{n}'$
are complex orientable cohomology theories carrying formal group laws
of particular interest.  The question of whether ${\rm BP}$ admits an
$E_\infty$ ring structure has been an open question for over 30 years,
and has become increasingly relevant in modern chromatic homotopy
theory.  We find that the connective cover of the cohomology theory
associated to an Atkin-Lehner quotient at discriminant $10$ is a model
for $\comp{(\BPP{2}')}_3$.  To our knowledge, this is the first known
example of a generalized $\BPP{n}$ having an $E_\infty$ ring structure
for $n > 1$.

At discriminant $10$, we obtain an example with torsion in its
homotopy groups.  We examine the homotopy of the spectrum
associated to the Shimura curve, and specifically at $p=3$ we give a
complete computation together with the unique lifts of the
Atkin-Lehner involutions.  The invariants under the Atkin-Lehner
involutions in the cohomology are isomorphic to the cohomology of the
$3$-local moduli of elliptic curves, and the associated homotopy fixed
point spectrum has the same homotopy groups as the localization
$\TMF_{(3)}$ of the spectrum of topological modular forms.  This
isomorphism does not appear to arise due to geometric considerations
and the reason for it is, to be brief, unclear.

The main technique allowing these computations to be carried out is
the theory of complex multiplication.  One first determines data for
the associated orbifold over $\mb C$, for example by simply taking
tables of this data \cite{alsinabayer}, or making use of the
Eichler-Selberg trace formula.  One can find defining equations for
the curve over $\mb Q$ using complex multiplication points and level
structures.  Our particular examples were worked out by Kurihara
\cite{kurihara} and Elkies \cite{elkies}.  We can then determine
uniformizing equations over $\mb Z[1/N]$ by excluding the possibility
that certain complex-multiplication points can have common reductions
at most primes.  This intersection theory on the associated arithmetic
surface was examined in detail by Kudla, Rapoport, and Yang \cite{kry}
and reduces in our case to computations of Hilbert symbols.  The rings
of automorphic forms can then be explicitly determined away from $2$
in terms of meromorphic sections of the cotangent bundle.  In order to
compute the Atkin-Lehner operators, we apply the Eichler-Selberg trace
formula.  Finally, we move from automorphic forms to homotopy theory
via the Adams-Novikov spectral sequence.

We now sketch the structure of the paper.  In
Section~\ref{sec:background} we review background material on
quaternion algebras, Shimura curves, Atkin-Lehner involutions, points
with complex multiplication, and automorphic forms.  In
Section~\ref{sec:spectra} we describe how Shimura curves give rise to
spectra via Lurie's theorem, and how their homotopy groups are related
to rings of automorphic forms.  In Section~\ref{sec:disc6} we find
defining equations for the Shimura curve of discriminant 6 over $\mb
Z[1/6]$ in detail by the sequence of steps described above, as well as
the ring of automorphic forms and the homotopy of various associated
spectra.  In Section~\ref{sec:disc14} we carry out the same program
for an Atkin-Lehner quotient of the Shimura curve of discriminant
$14$, culminating in the construction of an $E_\infty$-ring structure
on a generalized truncated Brown-Peterson spectrum.  In
Section~\ref{sec:disc10} we apply the same methods to the curve of
discriminant $10$.  The existence of elliptic points of order $3$, and
hence $3$-primary torsion information, complicates this computation
and necessitates working equivariantly with smooth covers induced by
level structures.  In Sections~\ref{sec:disc10spectra},
\ref{sec:disc10tate}, and \ref{sec:disc10fixedpoint} we use this
equivariant data to compute the homotopy groups of the resulting
spectrum via homotopy-fixed-point spectral sequences.

The authors would like to thank Mark Behrens, Benjamin Howard,
Niko Naumann, and John Voight for conversations related to this
material.

\section{Background}
\label{sec:background}

In this section we begin with a brief review of basic material on
Shimura curves.  Few proofs will be included, as much of this is
standard material in number theory.

For a prime $p$ and a prime power $q = p^k$, we write $\mb F_q$ for
the field with $q$ elements, $\mb Z_q$ for its ring of Witt vectors,
and $\mb Q_q$ for the associated field of fractions.  The elements
$\omega$ and $i$ are fixed primitive 3rd and 4th roots of unity in
$\mb C$.

\subsection{Quaternion algebras}

A quaternion algebra over $\mb Q$ is a $4$-dimensional simple algebra
$D$ over $\mb Q$ whose center is $\mb Q$.  For each prime $p$, the
tensor product $D \otimes \mb Q_p$ is either isomorphic to $M_2 \mb
Q_p$ or a unique division algebra
\begin{equation}
  \label{eq:divrel}
D_p = \{a + bS\ |\ a,b \in \mb Q_{p^2}\},
\end{equation}
with multiplication determined by the multiplication in $\mb Q_{p^2}$
and the relations $S^2 = -p$, $aS = S a^\sigma$.  (Here $a^\sigma$ is
the Galois conjugate of $a$.)  We refer to $p$ as either {\em split} or
{\em ramified} accordingly.  Similarly, the ring $D \otimes \mb R$ is
either isomorphic to $M_2 \mb R$ or the ring $\mb H$ of quaternions.
We correspondingly refer to $D$ as split at $\infty$ ({\em
indefinite}) or ramified at $\infty$ ({\em definite}).

We have a bilinear Hilbert symbol $(-,-)_p\co (\mb Q_p^\times)^2 \to
\{\pm 1\}$.  The element $(a,b)_p$ is $1$ or $-1$ according to whether
the algebra generated over $\mb Q_p$ by elements $i$ and $j$ with
relations $i^2 = a$, $j^2 = b$, $ij = -ji$ is ramified or not.

\begin{lem}[\cite{serrearith}]
Write $a = p^{\nu_p(a)} u$ and $b = p^{\nu_p(b)}v$.  If $p$ is odd,
we have
\[
(a,b)_p = \lsym{-1}{p}^{\nu_p(a) \nu_p(b)} \lsym{\bar u}{p}^{\nu_p(b)}
\lsym{\bar v}{p}^{\nu_p(a)},
\]
with $\bar x$ being the mod-$p$ reduction.  At $p=2$, we have
\[
(a,b)_2 = (-1)^{\frac{(u-1)(v-1)}{8} (2 + (u+1)\nu_2(a)  +  (v+1)\nu_2(b))}
\]
At $p = \infty$, we have $(a,b)_\infty$ equal to $-1$ if and only if
$a$ and $b$ are both negative.
\end{lem}
Here $\lsym{x}{p}$ is the Legendre symbol.

\begin{thm}[Quadratic reciprocity]
The set of ramified primes of a quaternion algebra $D$ over $\mb Q$ is
finite and has even cardinality.  For any such set of primes, there
exists a unique quaternion algebra of this type.
\end{thm}

We refer to the product of the finite, ramified primes as the
discriminant of $D$.  $D$ is a division algebra if and only of this
discriminant is not $1$.

The left action of any element $x$ on $D$ has a characteristic
polynomial of the form $(x^2 - \Tr(x) x +\Nm(x))^2$, and $x^2 - \Tr(x)
+ \Nm(x)$ is zero in $D$.  We define $\Tr(x)$ and $\Nm(x)$ to be the
{\em reduced trace} and {\em reduced norm} of $x$.  This trace and
norm are additive and multiplicative respectively.  There is a
canonical involution $x \mapsto x^\iota = \Tr(x) - x$ on $D$.

An {\em order} in $D$ is a subring of $D$ which is a lattice, and the
notion of a {\em maximal order} is clear.  The reduced trace and
reduced norm of an element $x$ in an order $\Lambda$ must be integers.
Two orders $\Lambda$, $\Lambda'$ are {\em conjugate} if there is an
element $x \in D$ such that $x\Lambda x^{-1} = \Lambda'$.

\begin{thm}[\cite{eichler}]
\label{thm:allconjugate}
If $D$ is indefinite, any two maximal orders $\Lambda$ and $\Lambda'$
are conjugate.
\end{thm}
Theorem~\ref{thm:allconjugate} generalizes the Noether-Skolem theorem
that any automorphism of $D$ is inner.  It is false for definite
quaternion algebras.

\begin{prop}
An order is maximal if and only if the completions $\Lambda
\otimes \mb Z_p$ are maximal orders of $D \otimes \mb Q_p$ for all
primes $p$.
\end{prop}
At split primes, this is equivalent to $\Lambda \otimes
\mb Z_p$ being some conjugate of $M_2(\mb Z_p)$.  At ramified primes,
there is a {\em unique} maximal order of the division algebra of
Equation~\ref{eq:divrel} given by
\begin{equation}
  \label{eq:orderrel}
\Lambda_p = \{ a + bS\ |\ a, b \in \mb Z_{p^2}\}.
\end{equation}

Let $F$ be a quadratic extension of $\mb Q$.  Then we say the field
$F$ {\em splits D} if any of the following equivalent statements hold:
\begin{itemize}
\item There is a subring of $D$ isomorphic to $F$.
\item $D \otimes F \cong M_2(F)$.
\item No prime ramified in $D$ is split in $F$.
\end{itemize}

\subsection{Shimura curves}

Fix an indefinite quaternion algebra $D$ over $\mb Q$ with
discriminant $N \neq 1$ and a maximal order $\Lambda \subset D$.  We can
choose an embedding $\tau\co D \into M_2(\mb R)$, and get a composite
map $\Lambda^\times \to \PGL_2(\mb R)$.  The latter group acts on the
space $\mb C \setminus \mb R$, which is the union ${\cal H} \cup
\overline{{\cal H}}$ of the upper and lower half-planes.  The
stabilizer of ${\cal H}$ is the subgroup $\Gamma = \Lambda^{\Nm=1}$ of
norm-$1$ elements in $\Lambda$, which is a cocompact Fuchsian
group.

\begin{defn}
The complex Shimura curve of discriminant $N$ is the the compact
Riemann surface ${\cal X}^D_{\mb C} = {\cal H} / \Gamma$.
\end{defn}

The terminology is, to some degree, inappropriate, as it should
properly be regarded as a stack or orbifold.  When the division
algebra is understood we simply write ${\cal X}_{\mb C}$.

This object has a moduli-theoretic interpretation.  The embedding
$\tau$ embeds $\Lambda$ as a lattice in a $4$-dimensional real vector
space $M_2(\mb R)$, and there is an associated quotient torus $\mb T$
with a natural map $\Lambda \to \End(\mb T)$.  The space ${\cal X}_{\mb
  C}$ is the quotient of the space of compatible complex structures by
the $\Lambda$-linear automorphism group, and parametrizes
$2$-dimensional abelian varieties with an action of $\Lambda$.  Such
objects are called {\em fake elliptic curves}.

The complex Shimura curve can be lifted to an algebraically defined
moduli object.  As in \cite{milne}, we consider the moduli of
$2$-dimensional projective abelian schemes $A$ equipped with an action
$\Lambda \to \End(A)$.  We impose a constraint on the
tangent space $T_e(A)$ at the identity: for any such $A/\Spec(R)$ and
a map $R \to R'$ such that $\Lambda \otimes R' \cong M_2(R')$, the two
summands of the $M_2(R')$-module $R' \otimes_R T_e(A)$ are locally
free of rank $1$.
\begin{thm}[{\cite[Section 14]{boutot}}]
  There exists a smooth, proper Deligne-Mumford stack ${\cal X}^D$
  over $\Spec(\mb Z[1/N])$ parametrizing $2$-dimensional abelian
  schemes with $\Lambda$-action.  The complex points form the curve
  ${\cal X}^D_{\mb C}$.
\end{thm}

\begin{rmk}
\label{rmk:polarization}
We note that related PEL moduli problems typically include the data of
an equivalence class of polarization $\lambda\co A \to A^\vee$
satisfying $x^\vee \lambda = \lambda x^\iota$ for $x \in \Lambda$.  In
the fake-elliptic case, there is a unique such polarization whose
kernel is prime to the discriminant, up to rescaling by $\mb Z[1/N]$.
\end{rmk}

\begin{thm}[{\cite{shimurareal}}]
The curve ${\cal X}$ has no real points.
\end{thm}

The object ${\cal X}$ has a model that is $1$-dimensional and proper
over $\Spec(\mb Z)$ and smooth over $\Spec(\mb Z[1/N])$, but has
singular fibers at primes dividing the discriminant.  We may view this
either as a curve over $\Spec(\mb Z[1/N])$ or an arithmetic surface.

\subsection{Atkin-Lehner involutions}
\label{sec:al}

An abelian scheme $A$ over a base scheme $S$ has a subgroup scheme of
$n$-torsion points $A[n]$; if $n$ is invertible in $S$, the map $A[n]
\to S$ is \'etale and the fibers over geometric points are isomorphic
to $(\mb Z/n)^{2\dim(A)}$.  We suppose that $A$ is $2$-dimensional and
has an action of $\Lambda$, so that we inherit an action of the ring
$\Lambda/n$ on $A[n]$.

Let $p$ be a prime dividing $N$.  There is a two-sided ideal $I$ of
$\Lambda$ such that $I^2 = (p)$; in the notation of
Equation~\ref{eq:orderrel}, it is the intersection of $\Lambda$ with
the ideal of $\Lambda_p$ generated by $S$.  The scheme $A[p]$ has a
subgroup scheme $H_p$ of $I$-torsion, and any choice of generator
$\pi$ of $I_{(p)}$ gives an exact sequence
\[
0 \to H_p \to A[p] \overto^\pi H_p \to 0.
\]
Hence the $I$-torsion is a subgroup scheme of rank $p^2$.

For $d | N$, the sum of these over $p$ dividing $d$ is the unique
$\Lambda$-invariant subgroup scheme $H_d \subset A$ of rank $d^2$, and
there is a natural quotient map $A \to A/H_d$.  It is of degree $d^2$,
$\Lambda$-linear, and these properties characterize $A/H_d$ uniquely
up to isomorphism under $A$.  This gives a natural transformation
\[
w_d\co A \mapsto (A/H_d)
\]
on this moduli of fake elliptic curves, and the composite $w_d^2$
sends $A$ to $A/A[d]$, which is canonically isomorphic to $A$ via the
multiplication-by-$d$ map $m_d\co A/A[d] \to A$.

\begin{prop}
The map $w_d$ gives rise to an involution on the moduli object
called an Atkin-Lehner involution. These involutions satisfy $w_d w_{d'}
\simeq w_{dd'}$ if $d$ and $d'$ are relatively prime.
\end{prop}

The Atkin-Lehner involutions form a group isomorphic to $(\mb Z/2)^k$,
with a basis given by involutions $w_p$ for primes $p$ dividing $N$.
There are associated quotient objects of ${\cal X}$ by any subgroup of
Atkin-Lehner involutions, and we write ${\cal X}^*$ for the full
quotient.

\begin{thm}[\cite{morita}]
For any group $G$ of Atkin-Lehner involutions, the coarse moduli
object underlying ${\cal X}/G$ is smooth.
\end{thm}

\begin{rmk}
\label{rmk:alnormalize}
The involution $w_d$ on ${\cal X}$ does not lift to an involution on the
universal abelian scheme ${\cal A}$, nor on the associated
$p$-divisible group ${\cal A}[p^\infty]$.  The natural isogeny $A \to
A/H_d$ of abelian varieties with $\Lambda$-action gives rise to a
natural $\Lambda$-linear isomorphism of $p$-divisible groups
$t'_d\co {\cal A}[p^\infty] \to (w_d)^* {\cal A}[p^\infty]$ over
${\cal X}$.  However, the natural diagram
\xym{
A \ar[drr]_{[d]} \ar[r]^{w_d} & A/H_d \ar[r]^{w_d} & A/A[d] \ar[d]^{m_d}_\sim\\
&&A
}
gives rise to a composite $\Lambda$-linear isomorphism of
$p$-divisible groups on ${\cal X}$ as follows:
\xym{
{\cal A}[p^\infty] \ar[r]^{t'_d}  \ar[drr]_{d} &
(w_d)^* {\cal A}[p^\infty] \ar[r]^{t'_d} &
(w_d^2)^*{\cal A}[p^\infty] \ar[d]^{m_d} \\
&&{\cal A}[p^\infty] 
}
This shows that under the canonical isomorphism of $(w_d^2)^* {\cal
A}$ with ${\cal A}$, the map $(t'_d)^2$ acts as the
multiplication-by-$d$ map on the $p$-divisible group rather than the 
identity.  Therefore, the universal $p$-divisible group does not
descend naturally to ${\cal X}/\langle w_d\rangle$.  However, if $x
\in \mb Z_p^\times$ satisfies $x^2 = \pm d$, we can lift the
involution $w_d$ of ${\cal X}$ to a $\Lambda$-linear involution $t_d =
x^{-1} t'_d$ of the $p$-divisible group and produce a $\Lambda$-linear
$p$-divisible group on the quotient stack.
\end{rmk}

\subsection{Complex multiplication in characteristic $0$}

Fix an algebraically closed field $k$.  A $k$-point of ${\cal X}$ consists
of a pair $(A,\phi)$ of an abelian surface over $\Spec(k)$ and an
action $\phi\co \Lambda \to \End(A)$.  We write $\End_\Lambda(A)$ for
the ring of $\Lambda$-linear endomorphisms.  The ring $\End(A)$ is
always a finitely generated torsion-free abelian group, and hence so
is $\End_\Lambda(A)$.

\begin{prop}[{\cite[pg. 202]{mumford}}]
  If $k$ has characteristic zero, then either $\End_\Lambda(A) \cong
  \mb Z$ or $\End_\Lambda(A) \cong {\cal O}$ for an order ${\cal O}$
  in a quadratic imaginary field $F$ that splits $D$.  In the latter
  case, we say $A$ has complex multiplication by $F$ (or ${\cal O}$).
\end{prop}

Complex multiplication points over $\mb C$ have explicit
constructions.  Fix a quadratic imaginary subfield $F$ of $\mb C$.
Suppose we have an embedding $F \to D^{op}$, and let ${\cal O}$ be the
order $F \cap \Lambda^{op}$.  The ring $\Lambda^{op}$ acts
$\Lambda$-linearly on the torus $\mb T = \Lambda \otimes S^1$ via
right multiplication.  There are precisely two (complex-conjugate)
complex structures compatible with this action, but only one
compatible with the right action of $F \otimes \mb R = \mb C$.  This
embedding therefore determines a unique point on the Shimura curve
with ${\cal O}$-multiplication.  The corresponding abelian surface $A$
splits up to isogeny as a product of two elliptic curves with complex
multiplication by $F$.

\begin{lem}
  Suppose $A$ represents a point on ${\cal X}_{\mb C}$ fixed by an
  Atkin-Lehner involution $w_d$.  Then $A$ has complex multiplication
  by $\mb Q(\sqrt{-d})$ if $d > 2$, and by $\mb Q(\sqrt{-1})$ or $\mb
  Q(\sqrt{-2})$ if $d=2$.
\end{lem}

\begin{proof}
  The statement $(w_d)^* A \cong A$ is equivalent to the existence of an
  $\Lambda$-linear isogeny $x\co A \to A$ whose characteristic
  polynomial is of the form $x^2 + tx + d$, because the kernel of such
  an isogeny must be $H_d$.  This element generates the quadratic
  field $\mb Q(\sqrt{t^2-4d})$.  This must be an imaginary quadratic
  field splitting $D$, and so must not split at primes dividing $d$.
  A case-by-case check gives the statement of the lemma.
\end{proof}

\subsection{Automorphic forms}
\label{sec:automorph}

Suppose $p$ does not divide $N$, so that we may fix an isomorphism
$\Lambda \otimes \mb Z_p \cong M_2(\mb Z_p)$; this is equivalent to
choosing a nontrivial idempotent $e \in \Lambda \otimes \mb Z_p$.  For
any fake elliptic curve $A$ with $\Lambda$-action the $p$-divisible
group $A[p^\infty]$ breaks into a direct sum $e \cdot A[p^\infty]
\oplus (1-e) \cdot A[p^\infty]$ of two canonically isomorphic
$p$-divisible groups of dimension $1$ and height $2$, and the formal
component of $A[p^\infty]$ carries a $1$-dimensional formal group law.
Similarly, $e \cdot \Lie(A)$ is a $1$-dimensional summand of the Lie
algebra of $A$.

The Shimura curve ${\cal X}$ over $\mb Z_p$ therefore has two
naturally defined line bundles.  The first is the cotangent bundle
$\kappa$ over $\Spec(\mb Z[1/N])$.  The second is the split summand
$\omega = e \cdot z^* \Omega_{A/{\cal X}}$ of the pullback of the
vertical cotangent bundle along the zero section $z\co {\cal X} \to
{\cal A}$ of the universal abelian surface.  The sections of this
bundle can be identified with invariant $1$-forms on a summand of the
formal group of ${\cal A}$.  The Kodaira-Spencer theory takes the
following form.

\begin{prop}
There are natural isomorphisms $\kappa \cong \omega^2
\cong z^*(\wedge^2 \Omega_{A/{\cal X}})$.
\end{prop}

\begin{proof}
This identification of the tangent space is based on the deformation
theory of abelian schemes.  See \cite{oort} for an introduction to
the basic theory.  For space reasons we will only sketch some details.

Given a point on the Shimura curve represented by $A$ over $k$, the
relative tangent space over $\mb Z[1/N]$ of the moduli of abelian
varieties at $A$ is in bijective correspondence with the set of
lifts of $A$ to $k[\epsilon]/(\epsilon^2)$.  This set of
deformations is isomorphic to the group
\[
H^1(A,TA) \cong H^1(A,{\cal O}) \otimes T_e(A) \cong T_e(A^\vee)
\otimes T_e(A).
\]
Here $TA$ is the relative tangent bundle of $A$, which is naturally
isomorphic to the tensor product of the trivial bundle and the tangent
space at the identity $T_e(A)$, and $A^\vee$ is the dual abelian
variety.  The polarization mentioned in Remark~\ref{rmk:polarization}
makes $H^1(A, {\cal O}) \cong T_e(A^\vee)$ isomorphic to $T_e(A)$,
with the right $\Lambda$-action induced by the canonical involution
$\iota$ on $\Lambda$.  The set of deformations that admit lifts of the
$\Lambda$-action are those elements equalizing the endomorphisms $1
\otimes x$ and $x^\vee \otimes 1$ of $T_e(A^\vee) \otimes T_e(A)$ for
all $x \in \Lambda$.  The cotangent space of ${\cal X}$ at $A$, which
is dual to this equalizer, therefore admits a description as a tensor
product:
\[
\Omega_{A^\vee/R} \otimes_\Lambda \Omega_{A/k} \cong (\Omega_{A/k})^t
\otimes_{M_2(Z_p)} \Omega_{A/k}.
\]
The tangent space of $A$ being free of rank $2$ implies that this is
isomorphic to $\wedge^2 \Omega_{A/k}$.
\end{proof}

We view $\kappa$ and $\omega^2$ as canonically identified.  An {\em
automorphic form} of weight $k$ on ${\cal X}$ is a section of
$\omega^{\otimes k}$, and one of even weight is a section of
$\kappa^{\otimes k/2}$.

The isomorphism $\kappa \cong \omega^2$ is only preserved by
isomorphisms of abelian schemes.  For instance, the natural map $A
\mapsto A/A[n]$ gives rise to the identity self-map on the moduli via
the isomorphism $m_n \co A/A[n] \to A$, and acts trivially on the
cotangent bundle, but it acts nontrivially on the vertical cotangent
bundle $\Omega_{A/{\cal X}}$ of the universal abelian scheme ${\cal A}$.

The canonical negation map $[-1]\co {\cal A} \to {\cal A}$ acts by
negation on $\omega$, and hence we have the following.
\begin{lem}
For $k$ odd, the cohomology groups
\[
H^i({\cal X}, \omega^{k}) \otimes \mb Z[1/2]
\]
are zero.
\end{lem}

\subsection{Application to topology}
\label{sec:spectra}

As in Section~\ref{sec:automorph}, we fix a nontrivial idempotent $e$
in $\Lambda \otimes \mb Z_p$.  Given a point $A$ of the moduli ${\cal
  X}$ over a $p$-complete ring, Serre-Tate theory shows that the
deformation theory of $A$ is the same as the deformation theory of $e
\cdot A[p^\infty]$; see \cite[Section 7.3]{taf} for the basic
argument.  We may then apply Lurie's theorem.

\begin{thm}
There is a lift of the structure sheaf ${\cal O}$ on the \'etale
site of $\comp{{\cal X}}_p$ to a sheaf ${\cal O}^{der}$ of locally
weakly even periodic $E_\infty$ ring spectra, equipped with an
isomorphism between the formal group data of the resulting spectrum
and the formal part of the $p$-divisible group.  Associated to an
\'etale map $U \to \comp{{\cal X}}_p$, the Adams-Novikov spectral
sequence for the homotopy of $\Gamma(U,{\cal O}^{der})$ takes the
form
\[
H^s(U, \omega^{\otimes t}) \rightarrow \pi_{2t-s} \Gamma(U,{\cal
  O}^{der}).
\]
\end{thm}

In particular, when $U = \comp{{\cal X}}_p$ we have an associated
global section object which might be described as a ``fake-elliptic''
cohomology theory.

\begin{defn}
We define the {\em $p$-complete fake topological modular forms
  spectrum} associated to the quaternion algebra $D$ to be the
$p$-complete $E_\infty$ ring spectrum
\[
\TAF^D_p = \Gamma(\comp{{\cal X}}_p,{\cal O}^{der}).
\]
We define the {\em complete fake topological modular forms spectrum}
associated to the quaternion algebra $D$ to be the
profinitely-complete $E_\infty$ ring spectrum
\[
\prod_{p \nmid N}\TAF^D_p.
\]
\end{defn}

If some subgroup of the Atkin-Lehner operators $w_d$ on the
$p$-divisible group are given choices of rescaling to involutions at
$p$ (see Remark~\ref{rmk:alnormalize}), we obtain an action of a
finite group $(\mb Z/2)^r$ on this spectrum by $E_\infty$ ring maps,
with associated homotopy fixed point objects equivalent to global
section objects on the quotient stack.

We may attempt to define ``$p$-local'' and ``global'' objects by
lifting the natural arithmetic squares
\xym{
\comp{{\cal X}}_p \otimes \mb Q\ar[r] \ar[d] & \comp{{\cal X}}_p \ar[d] &
\prod \comp{{\cal X}}_p \otimes \mb Q\ar[r] \ar[d] & \prod \comp{{\cal X}}_p \ar[d] \\
{\cal X} \otimes \mb Q \ar[r] & {\cal X}_{(p)} &
{\cal X} \otimes \mb Q \ar[r] & {\cal X}
}
to diagrams of $E_\infty$ ring spectra on global sections, perhaps
after a connective cover.  The spectrum associated to ${\cal X}
\otimes \mb Z_p$ is constructed by Lurie's theorem, whereas it may be
possible to construct the associated rational object by methods
of multiplicative rational homotopy theory.  We will discuss such a
global section object in a specific case in
Section~\ref{sec:disc6spectrum}.

\begin{rmk}
It seems conceivable that a direct construction analogous to that
for topological modular forms \cite{tmfconstr} might be possible,
bypassing the need to invoke Lurie's theorem.  However, this
construction is unlikely to be simpler than a construction of $\TMF$
with level structure.

Additionally, given more data about the structure of the modular curve
over $\mb Q$, it seems plausible that the methods of \cite[Section
9]{tmfconstr} might provide a direct lift of the structure sheaf of
${\cal X}$ to $E_\infty$ ring spectra.  However, as the rings of
rational automorphic forms have more complicated structure than those
in the modular case this requires further investigation.  
\end{rmk}

\section{Discriminant $6$}
\label{sec:disc6}

In this section, $D$ denotes the quaternion algebra of discriminant $N
= 6$ and ${\cal X} = {\cal X}^D$ the associated Shimura curve.  Our goal
in this section is to compute the cohomology
\[
H^s({\cal X}, \kappa^{\otimes t}),
\]
together with the action of the Atkin-Lehner involutions.  Write
${\cal X} \to X$ for the map to the underlying coarse moduli.

A fundamental domain for the action of the group of norm-$1$ elements
$\Gamma = \Lambda^{N=1}$ on the upper half-plane is pictured as follows.
(See \cite{alsinabayer}.)
\xynm{
0;/r8pc/:
a(150)="v1",
c+<-0.5pc,-0.5pc>*{v_1},
(-.366,.366)="v2",
c+<0pc,-0.5pc>*{v_2},
(0,0.2679)="v3",
c+<0.5pc,-0.5pc>*{v_3},
(.366,.366)="v4",
c+<0pc,-0.5pc>*{v_4},
a(30)="v5",
c+<0.5pc,-0.5pc>*{v_5},
(0,1)="v6",
c+<0.5pc,0.5pc>*{v_6},
(-1.2,0);(1.2,0)**@{.},
(-1,0),c-<0pc,0.5pc>*{-1},
(-0.5,0),c-<0pc,0.5pc>*{-0.5},
(0,0),c-<0pc,0.5pc>*{0},
(0.5,0),c-<0pc,0.5pc>*{0.5},
(1,0),c-<0pc,0.5pc>*{1},
(0,0);(0,1.2)**@{.},
"v5";p+a(120),**{},"v6",{\ellipse^{}},
"v1";p+a(60),**{},"v6",{\ellipse_{}},
"v2";p+a(130),**{},"v1",{\ellipse^{}},
"v2";p+a(10),**{},"v3",{\ellipse_{}},
"v4";p+a(170),**{},"v3",{\ellipse^{}},
"v4";p+a(50),**{},"v5",{\ellipse_{}},
}

The edges meeting at $v_2$ are identified, as are the edges meeting at
$v_4$ and those meeting at $v_6$.  The resulting surface has 2
elliptic points of order 2, coming from $v_6$ and $\{v_1,v_3,v_5\}$,
and 2 elliptic points of order 3, coming from $v_2$ and $v_4$.  The
underlying surface has genus 0 and hyperbolic volume
$\frac{2\pi}{3}$.

This fundamental domain leads to a presentation of the group $\Gamma$:
\[
\Gamma = \left\langle \gamma_{v_2}, \gamma_{v_4}, \gamma_{v_6} \ |\
  \gamma_{v_2}^3 = \gamma_{v_4}^3 = \gamma_{v_6}^2 =
  (\gamma_{v_2}^{-1} \gamma_{v_6} \gamma_{v_4})^2 = 1 \right\rangle
\]
Here $\gamma_{v}$ generates the stabilizer of the vertex $v$.  In
particular, $\Gamma$ has a quotient map to $\mb Z/6$ that sends
$\gamma_{v_2}$ and $\gamma_{v_4}$ to 2 and $\gamma_{v_6}$ to 3; this
is the maximal abelian quotient of $\Gamma$.  Let $K$ be the kernel of
this map.  This corresponds to a Galois cover ${\cal X}' \to {\cal X}$
of the Shimura curve ${\cal X}$ by a smooth curve, with Galois group
cyclic of order $6$.  The Riemann-Hurwitz formula implies that ${\cal
  X}'$ is of genus $2$.  This cover can be obtained by imposing level
structures at the primes $2$ and $3$.

\subsection{Points with complex multiplication}

We would like to now examine some specific points with complex
multiplication on this curve.  We need to first state some of
Shimura's main results on complex multiplication \cite{shimura} and
establish some notation.

Fix a quadratic imaginary field $F$ with ring of integers ${\cal
  O}_F$.  Let $I(F)$ be the symmetric monoidal category of fractional
ideals of ${\cal O}_F$, i.e. finitely generated ${\cal
  O}_F$-submodules of $F$ with monoidal structure given by product of
ideals.  The morphisms consist of multiplication by scalars of $F$.
We write $\Cl(F)$ for this {\em class groupoid} of isomorphism classes
in $I(F)$, with cardinality $\card{\Cl(F)}$ and mass $\mass{\Cl(F)} =
\card{\Cl(F)}/|{\cal O}_F^\times|$.  There is an associated Hilbert
class field $H_F$ which is a finite unramified extension of $F$ with
abelian Galois group canonically identified with $\Cl(F)$.

The Shimura curve ${\cal X}_{\mb C}$ has a finite set of points associated
to abelian surfaces with endomorphism ring ${\cal O}_F$.  The
following theorem (in different language) is due to Shimura.
\begin{thm}[{\cite[3.2, 3.5]{shimura}}]
\label{thm:cpxmult}
Let $A$ be a point with complex multiplication by ${\cal O}_F$.  Then
$A$ is defined over the Hilbert class field $H_F$.  If $S$ is the
groupoid of points with ${\cal O}_F$-multiplication that are isogenous
to $A$, then $S$ is a complete set of Galois conjugates of $A$ over
$F$, and there is a natural equivalence of $S$ with the groupoid
$I(F)$ of ideal classes.  This equivalence takes the action of the
class group $\Cl(F)$ to the Galois action.
\end{thm}

In the specific case of discriminant $6$, the curve ${\cal X}_{\mb Q}$
has two elliptic points of order $2$ that have complex multiplication
by the Gaussian integers $\mb Z[i]$ and are the unique fixed points by
the involution $w_2$.  They must be defined over the class field $\mb
Q(i)$, and since the curve ${\cal X}$ has no real points they must be
Galois conjugate over $\mb Q$.  The involutions $w_3$ and $w_6$
interchange them.  We denote these points by $P_1$ and $P_2$.

Similarly, there are two elliptic points of order $3$ with complex
multiplication by $\mb Z[\omega]$ that are Galois conjugate, defined
over $\mb Q(\omega)$, and have stabilizers $w_3$.  We denote these
points by $Q_1$ and $Q_2$.

Finally, there are two ordinary points fixed by the involution $w_6$
with complex multiplication by $\mb Z[\sqrt{-6}]$ that are Galois
conjugate.  If $K$ is the field of definition, then $K$ must be
quadratic imaginary and $K(\sqrt{-6})$ must be the Hilbert class field
$\mb Q(\sqrt{-6},\omega)$, which has only $3$ quadratic subextensions
generated by $\sqrt{-6}$, $\omega$, and $\sqrt{2}$.  Therefore, these
CM-points must be defined over $\mb Q(\omega)$.  We denote these
points by $R_1$ and $R_2$.

\subsection{Reduction mod $p$ and divisor intersection}
\label{sec:intersect}

An abelian surface $A$ in characteristic zero with complex
multiplication has reductions at finite primes.  In this section, we
examine when two points with complex multiplication can have a common
reduction at a prime $p$ not dividing the discriminant $N$.

Suppose $K$ is an extension field of $\mb Q$ with ring of integers
${\cal O}_K$.  As ${\cal X}$ is a proper Deligne-Mumford stack, any
map $\Spec(K) \to {\cal X}$ extends to a unique map $\Spec({\cal O}_K)
\to {\cal X}$, representing canonical isomorphism classes of
reductions mod-$p$ of an abelian surface over $K$.  The associated
image of $\Spec({\cal O}_K)$ can be viewed as a horizontal divisor on
the arithmetic surface ${\cal X}$ over ${{\cal O}_K}$.  If the point
has complex multiplication by a quadratic imaginary field, then so do
all points on this divisor.

We wish to obtain criteria to check when two such CM divisors can have
nonempty intersection over $\Spec(\mb Z[1/N])$.  Our treatment is
based on the intersection theory developed in \cite{kry}.  However,
our needs for the theory are more modest, as we only want to
exclude such intersections.

\begin{prop}[{\cite[5.2]{milne}}]
Suppose $k$ is an algebraically closed field of characteristic $p >
0$, $p \nmid N$.  The endomorphism ring of a point $(A,\phi)$
of the Shimura curve ${\cal X}$ over $k$ must be one of the following
three types.
\begin{itemize}
\item $\End_\Lambda(A) = \mb Z$.
\item $\End_\Lambda(A) = {\cal O}$ for an order in a quadratic
imaginary field $F$.
\item $\End_\Lambda(A)$ is a maximal order in $D'$, where $D'$
  is the unique ``switched'' quaternion algebra ramified precisely at
  $\infty$, the characteristic $p$, and the primes dividing $N$.
\end{itemize}
\end{prop}
We refer to the third type of surface as {\em supersingular}, and in
this case $A$ is isogenous to a product of two supersingular elliptic
curves.  Over a finite field, the endomorphism ring of
$(A,\phi)$ is always larger than $\mb Z$.

\begin{prop}
\label{prop:intersect}
Suppose $A$ and $A'$ are CM-points in characteristic zero, with
endomorphism rings generated by $x \in \End_\Lambda(A)$ and $y
\in \End_\Lambda(A')$, with $d_x$ and $d_y$ the discriminants of
the characteristic polynomials of $x$ and $y$.  Assume these generate
non-isomorphic field extensions of $\mb Q$.  If the divisors
associated to $A$ and $A'$ intersect at a point in characteristic $p$,
then there exists an integer $m$ with
\[
0 > \Delta \mathop{=}^{def} (2m + \Tr(x) \Tr(y))^2 - d_x d_y
\]
such that the Hilbert symbols
\[
(\Delta,d_y)_q = (d_x,\Delta)_q
\]
are nontrivial precisely when $q = \infty$ or $q\,|\,pN$.
\end{prop}

\begin{proof}
As the divisors associated to $A$ and $A'$ intersect, they must have
a common reduction $\bar A \cong \bar A'$ in characteristic $p$ with
endomorphism ring containing both $x$ and $y$.  This forces the
object $\bar A$ to be supersingular.

Let $\Lambda' = \End_\Lambda(\bar A)$ be the associated maximal order
in the division algebra $D'$.  The subset $L = \{u + vx + wy\, |
\,u,v,w\in\mb Z\}$ has a ternary quadratic form induced by the reduced
norm, taking integral values:
\begin{eqnarray*}
  (u + vx + wy) &\mapsto& -\Nm_{D'}(u + vx + wy)\\
&=& -u^2 -\Nm(x) v^2 - \Nm(y) w^2 - \Tr(x) uv - \Tr(y) uw + m \cdot vw
\end{eqnarray*}
Here $m \in \mb Z$ is an unknown integer.  On trace-zero elements, this
quadratic form takes an element to its square.  Let $W \subset L
\otimes \mb Q$ be the subset of trace-zero elements (having $2u =
-v\Tr(x) -w\Tr(y)$); the quadratic form restricts on $W$ to the form
\[
(v,w) \mapsto \frac{d_x}{4}\cdot v^2 + \frac{d_y}{4} \cdot w^2 + \left(\frac{\Tr(x)
    \Tr(y)}{2} + m\right) \cdot vw.
\]

This quadratic form determines the division algebra $D'$ as follows.
We can choose two elements $\alpha$ and $\beta$ in $W \subset D'$
diagonalizing the form:
\[
v' \alpha + w' \beta \mapsto d_1 v'^2 + d_2 w'^2.
\]
In particular, we can complete the square in such a way as to obtain
the diagonalized form with $(d_1, d_2) = (d_x/4,-\Delta/4d_x)$, which is
equivalent to $(d_x,\Delta)$, or symmetrically complete in the other
variable.

The elements $\alpha$ and $\beta$ then satisfy $\alpha^2 = d_1$,
$\beta^2 = d_2$, and $\alpha^2 + \beta^2 = (\alpha + \beta)^2$, so
they anticommute.  The Hilbert symbols $(d_1, d_2)_p$ must
then be nontrivial precisely at the places where $D'$ is ramified.
\end{proof}

As $m$ ranges over integer values, Proposition~\ref{prop:intersect}
leaves a finite number of cases that can be checked by hand.  For
example, if $D$ is a division algebra ramified at a pair of distinct
primes $2 < p < q$, then no divisor associated to a fixed point of the
involution $w_p$ intersects any divisor associated to a fixed point of
$w_q$.

We now apply this to the curve of discriminant $6$.  In the previous
section we constructed points $P_i$ over $\mb Q(i)$ with
endomorphism ring $\mb Z[i]$, $Q_i$ over $\mb Q(\omega)$ with endomorphism
ring $\mb Z[\omega]$, and $R_i$ over $\mb Q(\omega)$ with endomorphism
ring $\mb Z[\sqrt{-6}]$.  The discriminants of these rings are $-4$,
$-3$, and $-24$.

\begin{prop}
\label{prop:nonintersect6}
The divisors $P_1 + P_2$, $Q_1 + Q_2$, and $R_1 + R_2$ are pairwise
nonintersecting on ${\cal X}$.
\end{prop}

\begin{proof}
Proposition~\ref{prop:intersect} implies that the divisors associated
to $P_i$ and $Q_j$ can only intersect in characteristic $p$ if there
exists an integer $m$ such that the Hilbert symbols $(-4,4m^2 - 12)_q$
are nontrivial precisely for $q \in \{2,3,p,\infty\}$.
\begin{itemize}
\item When $m^2 = 0$, $(-4,-12)_q \neq 1$ iff $q \in \{3,\infty\}$.
\item When $m^2 = 1$, $(-4,-8)_q \neq 1$ iff $q \in \{2,\infty\}$.
\item When $|m| > 1$, $(-4,4m^2-12)_\infty = 1$.
\end{itemize}
Therefore, no pair $P_i$ and $Q_j$ can intersect over $\Spec(\Zos)$.

Similarly, the intersections $P_i \cap R_k$ are governed by Hilbert
symbols of the form $(-4, 4m^2-96)_q$.  In order for this to be ramified
at $3$, we must have the second term divisible by $3$, which reduces
to the cases $(-4,-96)_q$ and $(-4,-60)_q$.  These are both nontrivial only
at $3$ and $\infty$.

Finally, the intersections $Q_j \cap R_k$ are governed by Hilbert
symbols of the form $(-3,4m^2-72)_q$.  If $m$ is not divisible by $3$,
the Hilbert symbol at $3$ is $\lsym{(2m)^2}{3} = 1$.  We can then
reduce to the cases $(-3,-72)_q$ and $(-3,-36)_q$, neither of which
can be nontrivial at any finite primes greater than $3$.
\end{proof}

\subsection{Defining equations for the Shimura curve}
\label{sec:definingeq6}

In this section we obtain defining equations for the underlying coarse
moduli $X$ of the Shimura curve ${\cal X}$ over $\mb Z[1/6]$.  Our
approach is based on that of Kurihara \cite{kurihara}.  These
equations are well-known, at least over $\mb Q$, but we will
illustrate the method.

In general, a general smooth proper curve ${\cal C}$ of genus zero
over a Dedekind domain $R$ may not be isomorphic to $\mb P^1_R$ for
several reasons.

First, a general curve may fail to be geometrically connected, or
equivalently the ring of constant functions may be larger than $R$.
In our case, the curve ${\cal X}$ is connected in characteristic zero,
and hence connected over $\Zos$.

Second, a general curve may fail to have any points over the field of
fractions $K$.  A geometrically connected curve of genus zero over
$\Spec(K)$ is determined by a non-degenerate quadratic form in three
variables over $K$ that has a solution if and only if the curve is
isomorphic to $\mb P^1_K$.  For example, $X$ itself has no real
points, and so is not isomorphic to $\mb P^1_{\mb Q}$.

However, the Galois conjugate points $P_j \in X_{\mb Q}$ are
interchanged by $w_2$ and $w_6$, and hence have a common image $P$,
defined over $\mb Q$, in $X/w_2$ and $X/w_6$.  Similarly, the
Galois conjugate points $Q_j$ have a common image $Q$ defined over
$\mb Q$ in $X/w_3$ and $X/w_6$, and $R_j$ have a common image
$R$ in $X/w_2$ and $X/w_3$.  The curve $X^* =X/\langle
w_2,w_3\rangle$ has points $P$,$Q$, and $R$ over $\mb Q$.  Therefore,
these quotient curves are all isomorphic to $\mb P^1$ over $\mb Q$.

Finally, if a general curve ${\cal C}$ of genus zero has a
$\Spec(K)$-point, we can fix the associated divisor and call it
$\infty$.  Let $G = {\cal O} \rtimes {\cal O}^*$ be the group of
automorphisms of $\mb P^1$ preserving $\infty$, or equivalently the
group of automorphisms of $\mb A^1$.  The object ${\cal C}$ is then a
form of $\mb P^1$ classified by an element in the Zariski cohomology
group $H^1_{Zar}(\Spec(R),G)$.  In general we have that
$H^1_{Zar}(\Spec(R),{\cal O})$ is trivial, and so such a form is
classified by an element of $H^1_{Zar}(\Spec(R),{\cal O}^*)$, which is
isomorphic to the ideal class group of $R$.  In particular:

\begin{lem}
\label{lem:uniformize}
If $R$ has trivial class group, any geometrically connected curve of
genus zero over $\mb P^1_R$ having a $K$-point is isomorphic to $\mb
P^1_R$.

Given two such curves ${\cal C}$ and ${\cal C}'$, suppose $P_i$ and
$P_i'$ ($i \in \{1,2,3\}$) are $\Spec(K)$-points on ${\cal C}$ and
${\cal C}'$ whose associated $\Spec(R)$-points are nonintersecting
divisors.  There exists a unique isomorphism ${\cal C} \to {\cal C}'$
over $R$ taking $P_i$ to $P_i'$.
\end{lem}

At each residue field $k$ of $R$, the associated map ${\cal C}_k \to
{\cal C}'_k$ is the uniquely determined linear fractional
transformation moving the reduction of $P_i$ to that of $P_i'$.

As $\Zos$ has trivial class group, Lemma~\ref{lem:uniformize} implies
that the four curves $X/w_2$, $X/w_3$, $X/w_6$, and $X^*$ are
noncanonically isomorphic to $\mb P^1_{\Zos}$.  The divisors fixed by
the involutions $w_2$, $w_3$, and $w_6$ are nonintersecting over
$\Spec(\Zos)$ by Proposition~\ref{prop:nonintersect6}, allowing
explicit uniformizations to be constructed.

We can choose a defining coordinate $z\co X^* \to \mb P^1$ such
that $z(P) = 0$, $z(Q) = \infty$, $z(R) = 1$.

As the divisors $\{P_1,Q,R\}$ are nonintersecting on $X/w_3$ over
$\mb Z[1/6,\omega]$ and $\{\sqrt{-3},0,\infty\}$ are distinct on
$\mb P^1$ over $\mb Z[1/6,\omega]$, Lemma~\ref{lem:uniformize} implies
that there is a unique coordinate $x\co X/w_3 \to
\mb P^1$ defined over $\mb Z[1/6,\omega]$ such that $x(P_1) =
\sqrt{-3}$, $x(Q) = 0$, and $x(R) = \infty$.  The involutions $w_2$
and $w_6$ fix $Q$ and $R$, and hence must send $x$ to $-x$.
This is invariant under the Galois action and thus lifts to an
isomorphism of coarse moduli over $\Zos$.  We have $1 + (3/x^2) = z$.

Similarly, as $\{P,Q_1,R\}$ are nonintersecting on $X/w_2$ over
$\mb Z[1/6,i]$, there is a unique coordinate $y$ such that $y(P) = 0$,
$y(Q_1) = i$, and $y(R) = \infty$, invariant under the Galois action and
thus defined over $\Zos$.  We have $1 + 1/(y^2+1) = z$.  The
involutions $w_3$ and $w_6$ send $y$ to $-y$.

We then have
\[
1 + 3/x^2 = 1 - 1/(y^2 + 1)\text{, or }x^2 + 3y^2 + 3 = 0.
\]
The vanishing locus of $x$ consists of the points $Q_i$, and the
vanishing locus of $y$ consists of the points $P_i$.  The Atkin-Lehner
involutions act by $w_2 x = -x, w_2 y = y$ and $w_3 x = x, w_3 y =
-y$.

The coordinates $x$ and $y$ generate the pullback curve $X$, and we
thus have the description of $X$ as the smooth projective curve with
homogeneous equation
\[
X^2 + 3Y^2 + 3Z^2 = 0.
\]
Equations equivalent to this were given in \cite{kurihara} and
\cite{elkies}.

\subsection{The ring of automorphic forms}

In this section we use the defining equations for $X$ to compute
the ring of automorphic forms on ${\cal X}$.

The graded ring $\oplus H^0({\cal X}, \omega^{k})$ of automorphic forms
can be understood via the coarse moduli object away from primes
dividing the order of elliptic points.

\begin{lem}
\label{lem:pullback}
Suppose ${\cal X}$ has a finite number of elliptic points $p_i$ of order
$n_i$ relatively prime to $M$.  Then pullback of forms induces an isomorphism
\[
H^0\left(X, \kappa_{X}^{t}\left(\sum \left\lfloor
    t(n_i-1)/n_i\right\rfloor p_i\right)\right) \otimes \mb Z[1/M]
\to H^0({\cal X}, \kappa_{{\cal X}}^{t}) \otimes \mb Z[1/M].
\]
\end{lem}
In other words, a section of the $t$-fold power of the canonical
bundle $\kappa$ of ${\cal X}$ is equivalent to a section of $\kappa^t$ on
the coarse moduli that has poles only along the divisors $p_i$ of
degree less than or equal to $\frac{t(n_i-1)}{n_i}$.  This can be
established by pulling back a meromorphic $1$-form $\omega$ on $X$
to an equivariant $1$-form on ${\cal X}$, and examining the zeros and
poles.  Any other equivariant meromorphic section differs by
multiplication by a meromorphic function on $X$.

For discriminant $6$, then, the group of automorphic forms of weight
$k = 2t$ on ${\cal X}$ is identified with the set of meromorphic
sections of $\kappa^t$ on the curve $\{X^2 + 3Y^2 + 3Z^2=0\}$ in $\mb
P^2_\Zos$ that only have poles of order at most $\lfloor{t/2}\rfloor$ along
the divisor $X=0$ and at most $\lfloor{2t/3}\rfloor$ along $Y=0$.

In particular, we have the coordinates $x = X/Z$ and $y = Y/Z$.  The
form $dx$ has simple zeros along $Y=0$ and double poles along $Z=0$.
We have the following forms:
\[
  U = \frac{dx^3}{xy^5},\ V = \frac{dx^2}{xy^3},\ W = \frac{dx^6}{x^3y^{10}}
\]
These have weights $6$, $4$, and $12$ respectively.  They generate all
possible automorphic forms, and are subject only to the relation
\[
U^4 + 3V^6 + 3W^2 = 0.
\]

\subsection{Atkin-Lehner involutions and the Eichler-Selberg trace formula}
\label{sec:trace}

In this section we compute the action of the Atkin-Lehner involutions
on the ring of automorphic forms on ${\cal X}$.  The main tool in this is
the Eichler-Selberg trace formula, which is similar in application to
the Riemann-Roch formula or the Lefschetz fixed-point formula.

Associated to a union of double cosets $\Gamma \alpha \Gamma$ in
$\Lambda$, there is a Hecke operator on the ring of
automorphic forms.  In the particular case where this union consists of
the set of elements of reduced norm $\pm n$, the associated Hecke
operator is denoted by $T(n)$.  If $d$ divides $N$, then the
Atkin-Lehner operator $w_d$ is the Hecke operator $T(d)$.

We state the Eichler-Selberg trace formula in our case, as taken from
Miyake with minor corrections \cite[Section 6]{miyake}.  For our
purposes, we will not require level structure, nor terms involving
contributions from cusps that occur only in the modular case.

\begin{thm}[{\cite[Theorem 6.8.4]{miyake}}]
\label{thm:trace}
The trace of the Hecke operator $T(n)$ on the vector space of complex
automorphic forms of even weight $k \geq 2$ is given by
\begin{eqnarray*}
\Tr(T(n)) &=& \delta_{2,k} + \epsilon_n \frac{k-1}{12} n^{k/2 - 1} \prod_{p | N} (p-1)\\
 &&- \sum_t \frac{\alpha^{k-1} - \beta^{k-1}}{\alpha - \beta}
\prod_{p | N} \left(1 -
   \lsym{d}{p}\right)
 \sum_{\cal O} \mass{\Cl({\cal O})}
\end{eqnarray*}
Here $t$ ranges over integers such that $d = t^2 - 4n$ is negative,
$\alpha$ and $\beta$ are roots of $X^2 + tX + n$, and ${\cal O}$
ranges over all rings of integral elements ${\cal O} \supset \mb
Z[\alpha]$ of the field $\mb Q(\sqrt{d})$.  The element $\delta_{2,k}$
is $1$ if $k=2$ and $0$ if $k > 2$, and $\epsilon_n$ is $0$ unless $n$
is a square.
\end{thm}

\begin{rmk}
For the unique order ${\cal O} = \mb Z + f {\cal O}_F$ of index $f$
of a quadratic imaginary field $\mb Q(\sqrt d)$, the mass of the ring
class groupoid ${\Cl({\cal O})}$ is given by
\[
\mass{{\Cl({\cal O})}} = \mass{\Cl({\cal O}_F)} \cdot f \prod_{p | f} \left(1 -
  p^{-1}\lsym{d}{p}\right).
\]
\end{rmk}

\begin{rmk}
A word of caution about normalizations is in order.

There are multiple possible choices of normalization for the action of
the Hecke operators depending on perspective.  In the formula of
Theorem~\ref{thm:trace}, the Hecke operators $T(n)$ are normalized so
that an element $ A =
\begin{bmatrix}
  a & b \\ c & d
\end{bmatrix}
\in M_2(\mb R)$ acts on a form of weight $k$ by
\[
f(z) \mapsto \frac{\det(A)^{k-1}}{(cz+d)^k} f\left(\frac{az+b}{cz+d}\right).
\]
In particular, the double coset generated by a scalar $n \in \mb Z$ acts
trivially on forms of weight $2$ and rescales functions.  This
formula is computationally simpler and closely related to zeta
functions, but is not suitable for our purposes.  The normalization we
require is
\[
f(z) \mapsto {\left(\frac{\det(A)}{cz+d}\right)}^k f\left(\frac{az+b}{cz+d}\right),
\]
which is compatible with the description of automorphic forms as sections
of a tensor power of the line bundle of invariant differentials on the
universal abelian scheme.  In particular, this normalization makes
the Atkin-Lehner operators $w_d$ into ring homomorphisms.
\end{rmk}

As an immediate application of the trace formula, we find the
following.
\begin{cor}
If $n > 3$, the self-map $w_n$ on the group of automorphic forms of
even weight $k > 0$ has trace
\[
d \cdot \delta_{2,k} + (-n)^{k/2} \prod_{p|N} \left(1 -
  \lsym{-n}{p}\right) \sum_{{\cal O}} \mass{\Cl({\cal O})}.
\]
\end{cor}
The traces of $w_2$ and $w_3$ do not admit as compact a formulation.

The Eichler-Selberg trace formula on the curve of discriminant $6$
implies that the Atkin-Lehner involutions on the spaces of automorphic
forms have the following traces.

\begin{tabular}{c|ccccccccccccc}
  $\kappa^t$&0&1&2&3&4&5&6&7&8&9&10&11\\
\hline
  id       &1&0&1&1&1&1&3&1&3&3&3 &3\\
  $w_2$    &1&0&$-2^2$&$2^3$&$2^4$&$-2^5$&$-2^6$&$2^7$&$2^8$&$-2^9$&$-2^{10}$&$2^{11}$\\
  $w_3$    &1&0&$-3^2$&$-3^3$&$3^4$&$3^5$&$3^6$&$-3^7$&$-3^8$&$-3^9$&$3^{10}$&$3^{11}$\\
  $w_6$    &1&0&$6^2$&$-6^3$&$6^4$&$-6^5$&$6^6$&$-6^7$&$6^8$&$-6^9$&$6^{10}$&$-6^{11}$
\end{tabular}

At various primes, as in Remark~\ref{rmk:alnormalize}, we can extend
the universal $p$-divisible group to quotients of ${\cal X}$ and obtain
further invariant subrings.  However, these depend on choices.

\begin{exam}
If $2$ and $3$ have roots in $\mb Z_p$ (i.e. $p \equiv
\pm 1 \mod 24$), we can lift these involutions $w_d$ of ${\cal X}$ to
involutions $t_d$ on the ring of automorphic forms.  The rescaling
$t_2 = (\sqrt 2)^{-1} w_2$ on forms of weight $2t$ has trace equal to
the trace of $w_2$ divided by $(\sqrt 2)^{2t}$, and similarly for
$t_3$ and $t_6$.  These being ring homomorphisms, it suffices to
determine their effects on the generators.  In particular, when $t=2$
we find that the form $V$ is fixed only by $t_6$, when $t=3$ we find
that the form $U$ is fixed only by $t_2$, and when $t=6$, knowing that
$U^2$ is fixed and $V^3$ is fixed only by $t_6$ we find that $W$
is fixed only by $t_3$.  (The multiples of the form $W$ are the only
ones having only zeros at the points $R_i$ interchanged by the
Atkin-Lehner operators, and hence $W$ can only be rescaled.)

The ring of invariants under the resulting action of $\mb Z/2 \times
\mb Z/2$ is generated by $U^2$, $V^2$, and $UVW$, as the element $W^2$
becomes redundant.  These generators are subject only to the relation
\[
(U^2)^3V^2 + 3U^2V^8 + 3(UVW)^2 = 0.
\]
\end{exam}

\begin{exam}
As complementary examples, at $p=5$ we have square roots
of $\pm 6$, and hence two lifts of the involution $w_6$ to the ring
of automorphic forms.  The involution $(\sqrt 6)^{-1} w_6$ fixes $V$
and $W$, negating $U$, and leaves the invariant subring
\[
\mb Z_5[U^2,V,W]/(U^2)^2 + 3V^6 + 3W^2.
\]
The involution $(\sqrt{-6})^{-1} w_6$ fixes $U$ and $V$,
negating $W$, and leaves instead the subring
\[
\mb Z_5[U,V].
\]
\end{exam}

\subsection{Higher cohomology}

Since the orders of the elliptic points are units, 
the only higher cohomology groups are determined by a Serre duality
isomorphism
\[
H^1({\cal X}, \kappa^t) \cong H^0({\cal X}, \kappa^{1-t}).
\]
Explicitly, the open subsets ${\cal X} \setminus \{P_j\}$ and ${\cal
X} \setminus \{Q_j\}$ give an affine cover of ${\cal X}$ and produce a
Mayer-Vietoris exact sequence on cohomology.  Taking invariants under
the action of the Atkin-Lehner involutions, we can also obtain exact
sequences for the higher cohomology of quotient curves.

Letting $R$ denote the graded ring of automorphic forms
\[
R = \mb Z[1/6,U,V,W]/(U^4 + 3V^6 + 3W^2),
\]
we have a Mayer-Vietoris exact sequence as follows.
\[
0 \to R \to U^{-1} R \oplus V^{-1} R \to (UV)^{-1} R \to \oplus_t
H^1({\cal X}, \kappa^t) \to 0
\]
In particular, we find that the element $W(UV)^{-1}$ gives rise to a
class in $H^1({\cal X}, \kappa)$ such that multiplication induces 
duality between $H^0$ and $H^1$.

\subsection{The associated spectrum}
\label{sec:disc6spectrum}

For a prime $p$, the Adams-Novikov spectral sequence for the
associated spectrum $\TAF^D_p$ is concentrated in filtrations $0$
and $1$, which determine each other by duality.  The spectral sequence
therefore collapses at $E_2$.

\begin{thm} In positive degrees,
\[
\pi_* \TAF^D_p \cong \mb Z_p[U,V,W]/(U^4 + 3V^6 + 3W^2) \oplus \mb
Z_p\{D\},
\]
where $|U| = 12$, $|V| = 8$, $|W| = 24$, and $D$ is a Serre duality
class in degree $3$ such that multiplication induces a perfect pairing
\[
\pi_t \TAF^D_p \otimes \pi_{3-t} \TAF^D_p \to \mb Z_p.
\]
\end{thm}

We remark that we can construct a connective $E_\infty$ ring spectrum
${\taf}^D_p$ as the pullback of a diagram of Postnikov
sections.
\xym{
  {\taf}^D_p \ar[r] \ar[d] & \TAF^D_p[0..\infty] \ar[d] \\
  \mb S_{p}[0..3] \ar[r] & \TAF^D_p[0..3]
}
(The lower-left corner is an Eilenberg-Mac Lane spectrum, as $\mb
S_p$ has no nontrivial homotopy groups below degree $7$.  We denote
it this way to show that the lower map exists as a map of Postnikov
sections.)  The lower map has cofiber $\Sigma^3 \eilm{\mb Z_p}$
with third homotopy group generated by the Serre duality generator
from $H^1({\cal X},\kappa)$.  The homotopy of $\taf^D_p$ is then
simply a subring of $\pi_* \TAF^D_p$ mapping isomorphically to the
ring of automorphic forms on $\comp{{\cal X}}_p$.  Other than
aesthetic or computational reasons, there is currently no particular
reason to prefer this connective cover.

We now discuss the construction of an integral global section object.
We have a map of spectra $\prod_{p \nmid N} \taf^D_p \to
\left(\prod_{p \nmid N} \taf^D_p\right)_{\mb Q}$ that, on homotopy
groups, is the map from the ring of automorphic forms on $\comp{{\cal
X}}$ to the ring of automorphic forms on the rationalization
$\comp{{\cal X}} \otimes \mb Q$.

This latter ring contains the ring of automorphic forms on ${\cal
X}_{\mb Q}$.  There is a unique homotopy type $\taf^D_{\mb Q}$ of
$E_\infty$-ring spectrum with these homotopy groups.  It can be
defined as the following left-hand pushout diagram of free $E_\infty$
algebras over the Eilenberg-Mac Lane spectrum $\eilm{\mb Q}$, with the
homotopy groups given at right:
\xym{
\mb P_{\mb Q}[S^{48}] \ar[r] \ar[d]
& \mb P_{\mb Q}[S^{12} \vee S^8 \vee S^{24}] \ar[d] &
\mb Q[R] \ar[r] \ar[d] & 
\mb Q[U,V,W] \ar[d] \\
\mb S_{\mb Q} \ar[r] & \taf^D_{\mb Q} &
\mb Q \ar[r] & \mb Q[U,V,W]/(U^4 + 3 V^6 + 3 W^2)
}
For any $E_\infty$ algebra $T$, this homotopy pushout diagram gives
rise to a Mayer-Vietoris sequence of homotopy classes of maps of
$E_\infty$-algebras out to $T$, as follows: 
\[
\pi_{49}(T) \to [\taf^D_{\mb Q},T]_{E_\infty} \to \pi_{12}(T) \times
\pi_8(T) \times \pi_{24}(T) \to \pi_{48}(T).
\]
In particular, as $\pi_{49}(\taf^D_{\mb Q})$ and $\pi_{49}(\taf^D_p)$
are both zero, this shows any other $E_\infty$ algebra with this ring
of homotopy groups accepts a weak equivalence from $\taf^D_{\mb Q}$,
and $(\prod \taf^D_p)_{\mb Q}$ accepts a unique homotopy class of map
such that the image on homotopy is the subring of rational automorphic
forms.

We may then define $\taf^D$ as the homotopy pullback in the resulting
diagram:
\xym{
\taf^D \ar[r] \ar[d]&
\prod_{p \nmid N} \taf^D_p \ar[d]\\
\taf^D_{\mb Q} \ar[r] &
(\prod_{p \nmid N} \taf^D_p)_{\mb Q}
}
The spectrum $\taf^D$ is a connective $E_\infty$-ring spectrum whose
homotopy groups form the ring $\mb Z[1/6,U,V,W]/(U^4 + 3V^6 + 3W^2)$
of modular forms on ${\cal X}^D$.

From this, we can form a homotopy pullback diagram of localizations:
\xym{
\TAF^D \ar[r] \ar[d] & U^{-1} \taf^D \ar[d] \\
V^{-1} \taf^D \ar[r] & (UV)^{-1} \taf^D
}
The spectrum $\TAF^D$ has natural maps to $\TAF^D_p$, each of
which is a $p$-completion map.

\section{Discriminant $14$}
\label{sec:disc14}

In this section, $D$ denotes the quaternion algebra of discriminant
$14$ and ${\cal X} = {\cal X}^D$ the associated Shimura curve.  The
complex curve ${\cal X}_{\mb C}$ is of genus $1$ and has only two
elliptic points, both of order $2$.  The computations are very similar
to those in discriminant $6$.  We will explain some of the main
details with the aim of understanding the full Atkin-Lehner quotient
${\cal X}^*$ at the prime $3$.

The Atkin-Lehner operator $w_2$ has four fixed points.  There are two
Galois conjugate points with complex multiplication by $\mb Z[i]$
defined over $\mb Q(i)$, and there are also two Galois conjugate points
with complex multiplication by $\mb Z[\sqrt{-2}]$ defined over $\mb
Q(\sqrt{-2})$.

The operator $w_7$ has no fixed points, as $\mb Q(\sqrt{-7})$ does not
split this division algebra.

The operator $w_{14}$ has four fixed points, each having complex
multiplication by $\mb Z[\sqrt{-14}]$ and all Galois conjugate.

Straightforward application of Proposition~\ref{prop:intersect} shows
that the divisors associated to points with complex multiplication by
$i$, $\sqrt{-2}$, and $\sqrt{-14}$ cannot intersect.

Elkies \cite[5.1]{elkies} constructs a coordinate $t$ on the curve
${\cal X}^*$ over $\mb Q$ taking the points with complex
multiplication by $i$, $\sqrt{-2}$, and $\sqrt{-14}$ to $\infty$, $0$,
and the roots of $16t^2 + 13t + 8$ respectively.  These are distinct
away from the prime $2$ and hence determine an integral coordinate on
${\cal X}^*$.  The images of the points fixed by Atkin-Lehner
operators give $4$ elliptic points, one of order $4$ and three of
order $2$.

The Eichler-Selberg trace formula implies that the Atkin-Lehner
operators have the following traces on forms of even weight $k = 2t > 2$ on
${\cal X}$.
\begin{itemize}
\item The trace of the identity (the dimension) is $t$ if $t$ is even
  and $t-1$ if $t$ is odd.
\item The trace of $w_2$ is $(-1)^t2^{t+1}$ if $t \equiv 0,1 \mod 4$
  and $0$ otherwise.
\item The trace of $w_7$ is $0$.
\item The trace of $w_{14}$ is $(-1)^t 2 \cdot 14^t$.
\end{itemize}
When $k=2$, $w_2$ has trace $-2$, $w_7$ has trace $7$, and $w_{14}$
has trace $-14$.

We have $\lsym{-2}{3} = \lsym{7}{3} = \lsym{-14}{3} = 1$, so at $p=3$ we may
renormalize these operators to involutions $t_d$ whose fixed elements
are automorphic forms on the quotient curve ${\cal X}^*$.  These
involutions have the following traces:
\begin{itemize}
\item The trace of $t_2$ is $2$ if $t \equiv 0,1 \mod 4$
  and $0$ otherwise.
\item The trace of $t_7$ is $0$.
\item The trace of $t_{14}$ is $2$.
\end{itemize}
When $k=1$ or $2$, these operators all have trace $1$.  Elementary
character theory implies that the group of automorphic forms of weight
$2t$ on ${\cal X}^*$ has dimension $1 + \lfloor{t/4}\rfloor$, from which
the following theorem results.

\begin{thm}
This ring of automorphic forms on $\comp{({\cal X}^*)}_{(3)}$ has the form
\[
\mb Z_3[U,V].
\]
The generator $U$ is in weight $2$, having a simple zero at the
divisors with complex multiplication by $\mb Z[i]$, and the second
generator $V$ in weight $8$ is non-vanishing on these divisors.
\end{thm}

The homotopy spectral sequence associated to this quotient Shimura
curve at $3$ collapses at $E_2$, and the dual portion is concentrated
in degrees $-11$ and below.  The connective cover $\taf^D_3$ has
homotopy concentrated in even degrees and is therefore complex
orientable, via a ring spectrum map ${\rm BP} \to \taf^D_3$.  The
homotopy of $\taf^D_3$ is a polynomial algebra on the image of $v_1$ in
$\pi_4$ (a unit times $U$ that vanishes on the height $2$ locus) and
the image of $v_2$ in $\pi_{16}$ (a combination of $U^4$ and a unit
times $V$ that is non-vanishing on the height $2$ locus).  The
composite ring homomorphism
\[
\mb Z_3[v_1,v_2] \into {\rm BP}_* \to \pi_* \taf^D_3
\]
is therefore an isomorphism, and the ring $\pi_* \taf^D_3$ is obtained
from $\comp{({\rm BP}_*)}_3$ by killing a regular sequence $v_3',
v_4', \ldots$ of polynomial generators.

\begin{thm}
There is a generalized truncated Brown-Peterson spectrum $\BPP{2}'$
admitting an $E_\infty$ ring structure at the prime $3$.
\end{thm}

\begin{proof}
The spectrum $\taf^{D,*}_3 = \Gamma(\comp{({\cal X}^*)}_3, {\cal
O}^{der})$ is an $E_\infty$-ring spectrum with the homotopy type
of a $3$-complete generalized Brown-Peterson spectrum.  Similar to
Section~\ref{sec:disc6spectrum}, we construct an $E_\infty$ ring
$\taf^{D,*}_{\mb Q}$ as the free $\eilm{\mb Q}$ algebra $\mb P_{\mb
Q}(S^4 \vee S^{16})$ on generators in degree $4$ and $16$, and the
images of the generators of the rational ring of automorphic forms
on ${\cal X}^*$ determine a unique homotopy class of map
$\taf^{D,*}_{\mb Q} \to (\taf^{D,*}_3)_{\mb Q}$.  We then construct
a homotopy pullback square of $E_\infty$ ring spectra as follows:
\xym{
\taf^{D,*}_{(3)} \ar[r] \ar[d] & \taf^{D,*}_3 \ar[d]\\
\taf^{D,*}_{\mb Q}\ar[r] & (\taf^{D,*}_3)_{\mb Q}
}
We have $\pi_* \taf^{D,*}{(3)} \cong \mb Z_{(3)}[U,V]$.  The resulting
spectrum is therefore complex orientable via a $p$-typical
orientation, and the composite map
\[
\mb Z_{(3)} [v_1, v_2] \into {\rm BP}_* \to \pi_* \taf^{D,*}
\]
is an isomorphism.
\end{proof}

\section{Discriminant $10$}
\label{sec:disc10}

In this section, $D$ denotes the quaternion algebra of discriminant
$10$ and ${\cal X} = {\cal X}^D$ the associated Shimura curve.  The
fundamental domain appears as follows.
\xynm{
0;/r8pc/:
(-.6828,0.8363)="v1",
*\dir{*},
c+<-0.5pc,-0.5pc>*{v_1},
(-0.1172,0.1435)="v2",
*\dir{*},
c+<-0.8pc,0pc>*{v_2},
(-0.0201,0.0246)="v3",
*\dir{*},
c+<-0.8pc,0pc>*{v_3},
(0.0201,0.0246)="v4",
*\dir{*},
c+<0.8pc,0pc>*{v_4},
(0.1172,0.1435)="v5",
*\dir{*},
c+<0.8pc,0pc>*{v_5},
(0.6828,0.8363)="v6",
*\dir{*},
c+<0.5pc,-0.5pc>*{v_6},
(-0.9,0);(0.9,0)**@{.},
(-0.5,0)*{|},c-<0pc,0.8pc>*{-0.5},
(0,0),c-<0pc,0.8pc>*{0},
(0.5,0)*{|},c-<0pc,0.8pc>*{0.5},
(0,1)*{-},c-<0.8pc,0pc>*{i},
(0,0);(0,1.2)**@{.},
"v1";p+a(-20),**{},"v2",{\ellipse_{}},
"v2";p+a(-30),**{},"v3",{\ellipse_{}},
"v3";p+a(30),**{},"v4",{\ellipse_{}},
"v6";p+a(-160),**{},"v5",{\ellipse^{}},
"v5";p+a(-150),**{},"v4",{\ellipse^{}},
"v1";p+a(40),**{},"v6",{\ellipse_{}},
}

The edges meeting at $v_2$ are identified, as are the edges identified
at $v_5$.  The top and bottom edges are identified as well.  The genus
of the resulting curve is $0$.  There are four elliptic points of
order 3 and none of order 2.

\subsection{Points with complex multiplication}

The involutions $w_2$, $w_5$, and $w_{10}$ each have two fixed points
on the curve ${\cal X}$, having complex multiplication by $\mb
Q(\sqrt{-2})$, $\mb Q(\sqrt{-5})$, and $\mb Q(\sqrt{-10})$
respectively.  The associated Hilbert class fields are $\mb
Q(\sqrt{-2})$, $\mb Q(\sqrt{-5},i)$, and $\mb
Q(\sqrt{-10},\sqrt{-2})$.  The fixed points are therefore defined over
$\mb Q(\sqrt{-2})$, $\mb Q(i)$, and $\mb Q(\sqrt{-2})$ due to
Theorem~\ref{thm:cpxmult}.

The curve has four elliptic points of order $3$, with complex
multiplication by $\mb Z[\omega]$, that are acted on transitively by
the Atkin-Lehner involutions.  Their field of definition is $\mb
Q(\omega)$ and they are Galois conjugate in pairs.  Each divisor
intersects its Galois conjugate at the prime $3$.

Proposition~\ref{prop:intersect} implies that the divisors associated
to the fixed points of the Atkin-Lehner involutions are
nonintersecting.  The divisors associated to the elliptic points
cannot intersect the fixed divisors of $w_2$ or $w_5$, and can
intersect those of $w_{10}$ only at the prime $3$.  Therefore, the
involution $w_{10}$ must act as Galois conjugation on the elliptic
points.

\subsection{Level structures}

On the curve of discriminant $10$, the elliptic points of order $3$
create torsion phenomena and obstruct use of the previous methods for
computation.  In order to determine higher cohomology, it is necessary
to impose level structure to remove the $3$-primary automorphism
groups.

Write ${\cal A} \to {\cal X}$ for the universal abelian surface over the
Shimura curve ${\cal X}$.  Associated to an ideal $I$ of $\Lambda$
satisfying $M\Lambda \subset I \subset \Lambda$ we can consider the
closed subobjects ${\cal A}[I] \subset {\cal A}[M] \subset {\cal A}$ of
$I$-torsion points and $M$-torsion points.  Over $\Spec(\mb Z[1/NM])$,
the maps ${\cal A}[I] \to {\cal A}[M] \to {\cal X}$ are \'etale.  We have
a natural decomposition by annihilating ideal:
\[
{\cal A}[I] \cong \coprod_{\Lambda \supset J \supset I} {\cal A}\langle J\rangle
\]
In particular, if $I$ is maximal then ${\cal A}\langle I\rangle$ is
${\cal A}[I]$ minus the zero section.

We will write ${\cal X}_1(I) = {\cal A}\langle I\rangle$ as an arithmetic
surface over $\Spec(\mb Z[1/NM])$.  It is a modular curve
parametrizing abelian surfaces with a chosen primitive $I$-torsion
point.  If $I$ is a two-sided ideal, the cover ${\cal X}_1(I) \to {\cal X}$ is
Galois with Galois group $(\Lambda/I)^\times$.  For sufficiently small
$I$ (such that no automorphism of an abelian surface can preserve such
a point), the object ${\cal X}_1(I)$ is represented by a smooth curve.

In particular, for $d$ dividing $N$ there is an ideal $I_d$ which is
the kernel of the map $\Lambda \to \prod_{p | d} \mb F_{p^2}$.  This
ideal squares to $(d)$, and so we write ${\cal X}_1(\sqrt d)$ for the
associated cover over $\mb Z[1/N]$ with Galois group $\prod_{p | d}
\mb F_{p^2}^\times$.

For a general two-sided ideal $I$, the surface ${\cal X}_1(I)$ is
generally not geometrically connected, or equivalently its ring of
constant functions may be a finite extension of $\mb Z[1/NM]$.

\begin{thm}[{\cite[3.2]{shimura}}]
\label{thm:basefield}
  Let $(n) = I \cap \mb Z$.  The ring of constant functions on the
  curve ${\cal X}(I)_{\mb Q}$ is the {\em narrow class field} $\mb
  Q(\zeta_n)$ generated by a primitive $n$'th root of unity.
\end{thm}

In the case of discriminant $10$, to find a cover of the curve with no
order-$3$ elliptic points we impose level structure at the prime $2$.
In this case, the choice of level structure is equivalent to the
choice of a {\em full} level structure, so we will drop the subscript.
Let ${\cal X}(\sqrt 2) \to {\cal X}$ be the degree $3$ Galois cover
with Galois group $\mb F_4^\times$.  The cover is a smooth curve of
genus $2$, and the map on underlying curves is ramified precisely over
the four divisors with complex multiplication by $\mb Z[\omega]$.  The
ring of constant functions on ${\cal X}(\sqrt 2)_{\mb Q}$ is $\mb Q$
by Theorem~\ref{thm:basefield}.

Shimura's theorem on fields of definition of CM-points has a
generalization to include level structures.  Suppose $\mf a \subset
{\cal O}_F$ is an ideal, and let $I(F,\mf a)$ be the groupoid
parametrizing fractional ideals $I$ of $F$ together with a chosen
generator of $I/\mf a \cong {\cal O}_F/\mf a$.  Write $\Cl(F,\mf a)$
for the associated {\em ray class group} under tensor product, with
mass
\[
\mass{\Cl(F,\mf a)} = \mass{\Cl(F)}\cdot \mass{({\cal O}/\mf a)^\times}
\]
and cardinality
\[
\card{\Cl(F,\mf a)} = \mass{\Cl(F,\mf a)} \cdot |(1 + \mf a)^\times|.
\]
There is an associated ray class field $H_{F,\mf a}$, an
extension of $F$ generalizing the Hilbert class field, with Galois
group $\Cl(F,\mf a)$.

\begin{thm}[{\cite[3.2, 3.5]{shimura}}]
Suppose we have an embedding $F \to D^{op}$ with
corresponding abelian surface $A$ acted on by ${\cal O}_F$, together with
a two-sided ideal $J$ of $\Lambda$ and a level $J$ structure on $A$.
Let $(n) = J \cap \mb Z$ and $\mf a = J \cap {\cal O}_F$.  Then $A$ is
defined over the compositum $\mb Q(\zeta_n) \cdot H_{F,\mf a}$
of the ray class field with the cyclotomic field, and the set of
elements isogenous to $A$ is isomorphic to the ray class groupoid
$I(F,\mf a)$ as $\Cl(F,\mf a)$-sets.
\end{thm}

For the curve of discriminant $10$, the points with complex
multiplication by $\mb Z[\omega]$, $\mb Z[\sqrt{-2}]$, $\mb
Z[\sqrt{-5}]$, and $\mb Z[\sqrt{-10}]$ must have associated ideals
$(2)$, $(\sqrt{-2})$, $(2,1+\sqrt{-5})$, and $(2,\sqrt{-10})$
respectively, with ray class fields $\mb Q(\omega)$, $\mb
Q(\sqrt{-2})$, $\mb Q(\sqrt{-5},i)$, and $\mb
Q(\sqrt{-10},\sqrt{-2})$.  In particular, the order-$3$ elliptic
points and the fixed points of $w_2$ have all lifts defined over the
same base field.

\subsection{Defining equations for the curve}

Let $Q_j$ denote the $4$ elliptic points of ${\cal X}_{\mb Q}$, and
$P_j$, $P_j'$, $P_j''$ the points with complex multiplication by
$\sqrt{-2}$, $\sqrt{-5}$, and $\sqrt{-10}$ respectively.  These points
have rational images $Q$, $P$, $P'$, and $P''$ on ${\cal X}^*$.
Elkies has shown \cite[Section 4]{elkies} that there exists a
coordinate $t$ on the underlying curve $X^*$ taking the following
values $t(Q) = 0$, $t(P) = \infty$, $t(P') = 2$, and $t(P'') = 27$.
As the divisors associated to $P$, $P'$, and $P''$ are nonintersecting
over $\mb Z[1/10]$, this extends to a unique isomorphism $X^* \to
\mb P^1_{\mb Z[1/10]}$.

There is a coordinate $y$ on $X/w_5$ such that $y(P) = 0$, $y(P'') =
\infty$, and without loss of generality we can scale $y$ so that
$y(P_j') = \pm 5i$.  Similarly, there exists a coordinate $z$ on
$X/w_2$ such that $z(P') = 0$, $z(P'') = \infty$, and $z(P_j) = \pm
5 \sqrt{-2}$.  These coordinates satisfy $y^2 = t - 27$ and $z^2 = 4 -
2t$, and generate all functions on the curve.

We find, after rescaling, the coarse moduli $X$ can defined by the
homogeneous equation
\[
2Y^2 + Z^2 + 2W^2 = 0
\]
over $\mb Z[1/10]$ \cite{kurihara}.

\subsection{Defining equations for the smooth cover}

Any element of norm $5$ in $\Lambda$ has the trivial conjugation
action on the quotient ring $\mb F_4$, and hence the involution $w_5$
commutes with the group of deck transformations of ${\cal X}(\sqrt 2)$.
The quotient curve $X(\sqrt 2)/w_5$ has genus $0$.  The involutions
$w_2$ and $w_{10}$ anticommute with the group of deck transformations.

The points $P_j$ with complex multiplication by $\sqrt{-2}$ each
have three lifts to ${\cal X}(\sqrt 2)$ defined over $\mb Q(\sqrt{-2})$.
The degree $5$ isogeny $P_1 \to P_2$ canonically identifies the
$2$-torsion of $P_1$ with that of $P_2$, and so the involution $w_5$
lifts to an involution interchanging Galois conjugates.  On the
quotient curve, we therefore have that the rational point $P$ of
${\cal X}/w_5$ has three rational points as lifts.  As $P$ has no
automorphisms, the divisors of the associated lifts are
nonintersecting.

Therefore, there exists a coordinate $u$ on $X(\sqrt{-2})/w_5$ over
$\mb Z[1/10]$ such that $u$ takes the lifts of $P$ to the points $\pm
1$ and $\infty$.  The order-$3$ group of deck transformations must
therefore consist of the linear fractional transformations
\[
u \mapsto \frac{3+u}{1-u} \mapsto \frac{u-3}{u+1}.
\]
The fixed points, which are the lifts of the points $Q_i$, are the
points satisfying $u^2 + 3 = 0$. The degree-$3$ function $\frac{u^3 -
  9u}{u^2-1}$ is invariant under these deck transformations and hence
descends to a coordinate on $X/w_5$.  As it takes the value
$\infty$ on $P$ and $\pm 3\sqrt{-3}$ on $Q_i$, it must be one of the
functions $\pm y$.  By possibly replacing $u$ with $-u$ we may assume
\[
y = \frac{u^3-9u}{u^2-1}.
\]
The function field is therefore generated by the coordinates $u$ and
$z$.  Expressing $y$ in terms of $u$, we find that these satisfy the
relation
\[
z^2(u^2-1)^2 = -2(u^2+1)(u^2+2u+5)(u^2-2u+5).
\]

\subsection{Automorphic forms}
\label{sec:autform10}

An automorphic form of weight $2t$ on ${\cal X}(\sqrt 2)$ is simply a
holomorphic section of $\kappa^t$ on $X(\sqrt 2)$.  Any of these is of
the form $f(u,z) du^t$, and to find generators one needs determine
when such a function is holomorphic.  The form $du$ itself has simple
zeros along the six divisors with $z = 0$, and double poles along the
two divisors with $u = \infty$.

The smooth cover ${\cal X}(\sqrt{2})$ has automorphic forms of weight $2$
given by
\[
U = \frac{u\,du}{z(u^2-1)}, V = \frac{du}{z(u^2 - 1)},
\]
together with an automorphic form of weight $6$ given by
\[
W = \frac{du^3}{(u^2+1)(u^2+2u+5)(u^2-2u+5)}.
\]
These generate the entire ring of forms.  These satisfy the formula
\[
W^2 = -2(U^2 + V^2)(U^2 + 2UV + 5V^2)(U^2 - 2UV + 5V^2),
\]
and the action of the deck transformation $\sigma\co u \mapsto
\frac{3+u}{1-u}$ is given as follows.
\begin{eqnarray*}
U &\mapsto& (-U-3V)/2\\
V &\mapsto& (U-V)/2\\
W &\mapsto& W
\end{eqnarray*}
We can then alternatively express the ring of automorphic forms on
this smooth cover as having generators $\{A,B,C,W\} = \{2V,2 \sigma V,
2 \sigma^2 V,W\}$, subject only to the relations $A+B+C=0$ and $W^2 =
-(A^2 + B^2)(B^2 + C^2)(C^2 + A^2)$.

We pause here to note the structure of the supersingular locus at $3$.
There are $2$ supersingular points at $3$ corresponding to the prime
ideal generated by $3$ and $U$, and $V$ vanishes at neither point.
Inverting $V$ and completing with respect to the ideal $(3,U)$, we
obtain the complete local ring
\[
\mb Z_3\pow{u}[w,V^{\pm 1}]/(w^2 + 2(u^2+1)(u^2 + 2u + 5)(u^2 - 2u + 5)).
\]
The generators are $u = U/V$ and $w = W/V^3$.  The power series $w^2$
has a constant term which is, $3$-adically, a unit and a square, and
the binomial expansion allows us to extract a square root.  Therefore,
this ring is isomorphic to
\[
\mb Z_3\pow{u}[V^{\pm 1}] \times \mb Z_3\pow{u}[V^{\pm 1}].
\]
This is a product of two summands of the Lubin-Tate ring, each
classifying deformations of one of the supersingular points.

\subsection{Atkin-Lehner involutions}

At the prime $3$, both $-2$ and $-5$ have square roots.  Therefore,
the $p$-divisible group on ${\cal X}$ can be made to descend to the curves
${\cal X}^*$ and ${\cal X}(\sqrt 2)/w_5$.  We write $t_2$ and $t_5$ for the
associated involutions on the rings of automorphic forms.

The involution $t_5$ negates $w = W/V^3$ and fixes $u = U/V$, and also
commutes with the action of the cyclic group of order $3$.  The lift
$t_5$ must either fix all forms of weight $2$ on ${\cal X}(\sqrt 2)$ and
negate $W$, or negate all forms of weight $2$ and fix $W$.  The former
must be the case: in the cohomology of ${\cal X}(\sqrt 2)/w_5$ the element
$U$ vanishing on the supersingular locus must survive to be a Hasse
form $v_1$ in weight $2$, as the element $\alpha_1$ in the cohomology
of the moduli of formal group laws has image zero.

Therefore, at $3$ the curve ${\cal X}(\sqrt 2)/t_5$ has
\[
\mb Z_3[U,V]
\]
as its ring of automorphic forms.

The involution $t_2$ conjugate-commutes with the group of deck
transformations, and choosing a lift gives rise to an action of the
symmetric group $\Sigma_3$ on this ring.  The fact that $t_2$ is
nontrivial forces the representation on forms of weight $2$ to be a
$\mb Z[1/10]$-lattice in the irreducible $2$-dimensional
representation of $\Sigma_3$.  There are two isomorphism classes of
such lattices.  Both have a single generator over the group ring $\mb
Z[1/10,\Sigma_3]$, and are distinguished by whether their reduction
mod $3$ is generated by an element fixed by $t_2$ or negated by
$t_2$.  In the former case, $H^1({\cal X}^*,\kappa) = 0$, whereas in
the latter it is $\mb Z/3$.  We will show in the following section
that a nonzero element $\alpha_1 \in H^1({\cal X}^*,\kappa)$ exists
due to topological considerations.

Therefore, the ring of $t_5$-invariant automorphic forms can be
expressed as a polynomial algebra
\[
\mb Z_3[A,B,C]/(A+B+C)
\]
on generators cyclically permuted by the group of deck transformations,
and such that one lift of $t_2$ sends $(A,B,C)$ to $(-A,-C,-B)$.

The ring of forms on ${\cal X}(\sqrt 2)$ invariant under the action of
$\mb Z/3$ is the algebra
\[
\mb Z_3[\dfour, \dsix, \sqrtdelta, \U]/(\U^2+2\dfour^3+\dsix^2,\, \sqrtdelta^2+4\dfour^3+27\dsix^2).
\]
Here $\dfour$, $\dsix$, and $\sqrtdelta$ are standard invariant
polynomials arising in the symmetric algebra on the reduced regular
representation generated by $A$, $B$, and $C$.  The subring of
elements invariant under the involutions $t_2$ and $t_5$ is
\[
\mb Z_3[\dfour, \sqrtdelta, \dsix^2]/(\sqrtdelta^2 + 4\dfour^3 + 27\dsix^2).
\]
Abstractly, this ring is isomorphic to the ring of integral modular
forms on the moduli of elliptic curves, with $\dsix^2$ playing the
role of the discriminant.

\subsection{Associated spectra}
\label{sec:disc10spectra}

The $\mb Z/3$-Galois cover of ${\cal X}$ by ${\cal X}(\sqrt 2)$ gives rise to
an expression of $\TAF^D_p$ as a homotopy fixed point spectrum of an
object $\TAF^D_p(\sqrt 2)$ by an action of $\mb Z/3$.  We will use the
homotopy fixed point spectral sequence to compute the homotopy groups
of $\TAF^D_p$. The difference between this spectral sequence
and the Adams-Novikov spectral sequence is negligible: the filtration
of the portion arising via Serre duality is lowered by $1$ in the
homotopy fixed point spectral sequence.

The positive-degree homotopy groups of $\TAF^D_p(\sqrt 2)$ form the
ring generated by the elements $A$, $B$, $C$, and $W$ of Section
\ref{sec:autform10}, together with a Serre duality class $D$ in
degree $3$ that makes the multiplication pairing perfect.

The main computational tool will be the Tate spectral sequence. The
nilpotence theorem \cite{nilpotence} guarantees that the Tate spectrum
is contractible (the periodicity class $\beta$ is
nilpotent). Therefore, the Tate spectral sequence, formed by inverting
the periodicity class, converges to $0$.  The Tate spectral sequence
has the added advantage that the $E_2$-term can be written in terms of
permanent cycles and a single non-permanent cycle.

\begin{thm}
The cohomology of $\mb Z/3$ with coefficients in the ring of
automorphic forms on ${\cal X}(\sqrt 2)$ is given by
\begin{multline*}
E_2=\Z_{(3)}[\dfour, \dsix, \U, \sqrtdelta][\beta]\otimes E(\alpha_1) / \\ (3,\dfour,\sqrtdelta)\cdot(\beta, \alpha_1),\, \U^2+2\dfour^3+\dsix^2,\, \sqrtdelta^2+4\dfour^3+27\dsix^2.
\end{multline*}
Here $\beta$ is the generator of $H^2(\Z/3;\Z)$ and $\alpha_1$ is
the element of $H^1$ coming from the reduced regular representation in
degree $4$.
\end{thm}

The Atkin-Lehner operators act as follows.  The involution $t_5$ fixes
$\alpha_1$, $\beta$, $\dfour$, $\dsix$, and $\sqrtdelta$, but negates
$\U$.  The involution $t_2$ fixes $\alpha_1$, $\dfour$, and
$\sqrtdelta$, but negates $\beta$ and $\dsix$.

This cohomology, $3$-completed, forms the portion of the Adams-Novikov
$E_2$-term for $\TAF^D_3$ generated by non-Serre dual classes.

\begin{cor}
For degree reasons, the first possible differentials are $d_5$ and $d_9$.
\end{cor}

The class $\alpha_1$ is so named because it is the Hurewicz image of
the generator of $\pi_3(S^0)$ at $3$.  This can be seen by
$K(2)$-localization: the $K(2)$-localization of $\TAF^D_3$ is a
wedge of four copies of $\EO_2(\Z/3\times \Z/2)$, corresponding to the
four supersingular points on ${\cal X}$. The image of $\alpha_1$ under the
composite
\[
\pi_3(S^0)\to\pi_3\TAF^D_3\to\pi_3\big(\EO_2(\Z/3\times\Z/2)\big)
\]
is the image of $\alpha_1$ under the Hurewicz homomorphism for
$\EO_2(\Z/3\times\Z/2)$, which is known
to be non-zero. By considering the map of Adams-Novikov spectral
sequences induced by the unit map $S^0\to\TAF^D_3$ and the
$K(2)$-localization map $\TAF^D_3\to L_{K(2)}\TAF^D_3$, we see that
the Adams-Novikov filtration of the element representing $\alpha_1$
can be at most $1$, allowing us to explicitly determine the element.
We note that (as required in the previous section) this also forces it
to be fixed by the Atkin-Lehner involutions.

\begin{prop}
The element $\beta_1=\langle\alpha_1,\alpha_1,\alpha_1\rangle$ is given by $\beta\dsix$.
\end{prop}
\begin{proof}
  Given a torsion-free commutative ring $R$ with $G$-action, a
  $1$-cocycle $f\co G \to R$ has a family of associated $1$-cochains
  $f^k\co g \mapsto f(g)^k$.  These satisfy $\delta f^2 = -2 (f
  \smallsmile f)$ and $\delta f^3 = -3\langle f,-2f,f\rangle$.  In
  particular, $f^3$ represents a cocycle mod $3$ whose image under the
  Bockstein map is in $\langle f,2f,f\rangle$.

  In the case in question, the $1$-cocycle $\alpha_1 \in H^1(\mb
  Z/3; \pi_4 \TAF^D_p(\sqrt 2))$ is represented by the trace-zero
  element $A$, and the element $A^3$ represents, in $H^1(Z/3; \pi_4
  \TAF^D_p(\sqrt 2)/3)$, the cube of the associated cocycle.  The
  image under the Bockstein is $\beta (A^3 + B^3 + C^3)/3 = -\beta
  \sigma_6$.  We then find $\beta \sigma_6 = \langle \alpha_1,
  \alpha_1, \alpha_1\rangle$ as desired.
\end{proof}

\begin{rmk}
  We can also detect $\beta_1$ using the $K(2)$-localization. There is
  only one class in the appropriate degree that is invariant under the
  Atkin-Lehner involutions, namely $\beta\dsix$. By an argument
  exactly like that for $\alpha_1$, we conclude that this is
  $\beta_1$.
\end{rmk}

\begin{cor}
The class $\beta\dsix$ is a permanent cycle.
\end{cor}

Most classes on the zero-line in the Adams-Novikov spectral sequence
are annihilated by multiplication by either $\alpha_1$ or $\beta$ ---
specifically, those classes that arise as elements in the image of the
transfer.  These play little role in the spectral sequence. However,
there are classes which play important roles for the higher
cohomology.

\begin{lem}
On classes of non-zero cohomological degree, multiplication by
$\dsix$ is invertible.
\end{lem}

\begin{proof}
  This follows easily from the cover of the moduli stack by affine
  spaces in which one of $\dfour$ or $\dsix$ is inverted. The
  Mayer-Vietoris sequence, together with the fact that
  $\alpha_1$ and $\beta$ are annihilated by $\dfour$, shows
  this result immediately.
\end{proof}

\begin{cor}
  We can choose a different polynomial generator for the higher
  cohomology: $\beta_1=\beta\dsix$.
\end{cor}

In particular, even though there is no class $\dsix^{-1}$, on classes
of non-zero cohomological degree the use of such notation is
unambiguous.

\begin{cor}
The Tate $E_2$-term is given by
\[
E_2=\F_9[\beta_1^{\pm 1},\dsix^{\pm 1}]\otimes E(\alpha_1).
\]
\end{cor}

Here the element $\sqrt{-1} \in \mb F_9$ is represented by the class
$\U/\dsix$, which must be a permanent cycle.  The involution $t_5$
fixes $\alpha_1$, $\beta_1$, and $\dsix$, but acts by conjugation on
$\mb F_9$.  The involution $t_2$ fixes $\alpha_1$ and $\beta_1$, but
negates $\dsix$ and acts by conjugation on $\mb F_9$.

\subsection{The Tate spectral sequence}
\label{sec:disc10tate}
\begin{figure}[h]
\includegraphics[angle=-90, width=\textwidth]{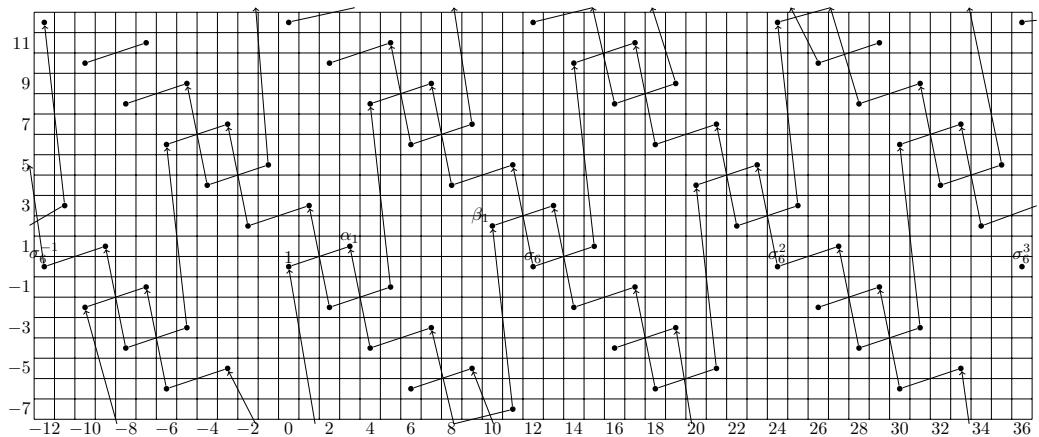}
\caption{The Tate $E_2$-term}
\label{fig:Tate}
\end{figure}

We suggest that the reader follow the arguments of this section with
the aid of Figure~\ref{fig:Tate}. In this picture, dots indicate a
copy of $\F_9$, lines of slope $1/3$ correspond to
$\alpha_1$-multiplication, and the elements $\dsix^i$ and $\beta_1$
are labeled. Both the $d_5$ and $d_9$-differentials are included,
allowing the reader to readily see that these are essentially the only
possibilities.

Though we use permanent cycles for our description of the higher
cohomology, it is more straightforward to determine differentials
originating from the periodicity class $\beta$.
\begin{thm}
There is a $d_5$-differential
\[
d_5(\beta)=\alpha_1\beta^3=\alpha_1\beta_1^3\dsix^{-3}.
\]
\end{thm}

\begin{proof}
  This is a relatively standard homotopy fixed point argument. The
  unit map $S^0\to \TAF^D_3$ is a $\Z/3$-equivariant map. This therefore
  induces a map of filtered spectra on the homotopy fixed point
  objects:
\[
D(B\Z/3_+)=(S^0)^{h\Z/3}\to \TAF^D_3
\]
This in turn gives rise to a map of homotopy fixed point spectral
sequences, where the source spectral sequence coincides with the
Atiyah-Hirzebruch spectral sequence for stable cohomotopy of
$B\Z/3$. We can therefore identify classes in the source spectral
sequence with homotopy classes and null-homotopies in the target
spectral sequence.

If we denote by $b_k$ the generator in $H^{2k}(\mb Z/3,\pi_0 S^0)$ on
the $E_2$-page for the source spectral sequence, then by definition
the image of $b_k$ in the homotopy fixed point spectral sequence for
$\TAF^D_3$ is $\beta^k$. Consideration of the cohomology of $B\mb
Z/3$ shows that in the Atiyah-Hirzebruch spectral sequence $b_1$
supports a differential to $\alpha_1$ times $b_3$.  We conclude that
in the homotopy fixed point spectral sequence for $\TAF^D_3$, the
class $\beta$ supports a differential whose target is
$\alpha_1\beta^3$.
\end{proof}

\begin{cor}
There is a $d_5$-differential of the form
\[
d_5(\dsix)=-\alpha_1\beta_1^2\dsix^{-1}.
\]
\end{cor}

\begin{proof}
  Since $\beta\dsix$ is a permanent cycle, we must have this
  differential by the Leibniz rule.
\end{proof}

\begin{rmk}
  An analogous differential appears in the computation of the homotopy of
  $\tmf$ \cite{bauer}, and we can repeat the argument here as
  well. The Toda relation $\alpha_1\beta_1^3=0$, since it arises from
  the sphere, must be visible here. For degree reasons, the only
  possible way to achieve this is to have a differential on $\dsix^2$
  hitting $\alpha_1\beta_1^2$, and this implies the aforementioned
  differential.
\end{rmk}

\begin{rmk}
  The differential on $\dsix$ also follows from geometric methods
  similar to those employed for $\beta$. In this case, instead of the
  unit map $S^0\to \TAF^D_3$, we use a multiplicative norm map $S^{2\rho}\to
  \TAF^D_3(\sqrt 2)$, where $\rho$ is the complex regular
  representation. The source homotopy fixed point spectral sequence
  is then the Atiyah-Hirzebruch spectral sequence for the
  Spanier-Whitehead dual of the Thom spectrum for a bundle over $B\mb
  Z/3$.
\end{rmk}

Since $\beta_1$ is an invertible permanent cycle, this leaves a
relatively simple $E_6$-page.

\begin{prop}
As an algebra,
\[
E_6=\F_9[\dsix^{\pm 3}, \beta_1^{\pm 1}]\otimes E([\alpha_1\dsix]).
\]
\end{prop}

A Toda style argument gives the remaining differential \cite{Greenbook}.
\begin{prop}
There is a $d_9$-differential of the form
\[
d_9(\alpha_1\dsix^4)=\toda{\alpha_1,\alpha_1\beta_1^2,\alpha_1\beta_1^2}=\beta_1^5.
\]
\end{prop}

Since $\beta_1$ is a unit, this gives a differential
\[
d_9(\alpha_1\dsix^4\beta_1^{-5})=1,
\]
so $E_{10}=0$.

\subsection{The homotopy fixed point spectral sequence}
\label{sec:disc10fixedpoint}
Having completed the Tate spectral sequence, the homotopy fixed
point spectral sequence and the computation of the homotopy of
$\TAF^D_3$ become much simpler. There is a natural map 
of spectral sequences from the homotopy fixed point spectral sequence to
the Tate spectral sequence that controls much of the behavior.  We
first analyze the kernel of the map
\[
H^\ast(\Z/3;H^*({\cal X}(\sqrt 2)))\to H^\ast_{Tate}(\Z/3;H^*({\cal
  X}(\sqrt 2))).
\]

\begin{lem}
  Everything in the ideal $(3, \dfour, \sqrtdelta)$ is a permanent cycle.
\end{lem}
\begin{proof}
  These classes are all in the image of the transfer.
\end{proof}
This shows that the differentials and extensions take place in
portions seen by the Tate spectral sequence, and we will henceforth work
in the quotient of the homotopy fixed point spectral sequence by this
ideal. By naturality, the canonical map of spectral sequences is an
inclusion through $E_5$, and the $d_5$-differential in the homotopy
fixed point spectral sequence is the $d_5$-differential in the Tate
spectral sequence.
\begin{cor}
\[
d_5(\dsix)=-\alpha_1\beta^2\dsix,\, d_5(\beta)=\alpha_1\beta^3,\text{ and }d_5(\U)=-\alpha_1\beta^2\U.
\]
\end{cor}

\begin{figure}[h]
\includegraphics[angle=-90, width=\textwidth]{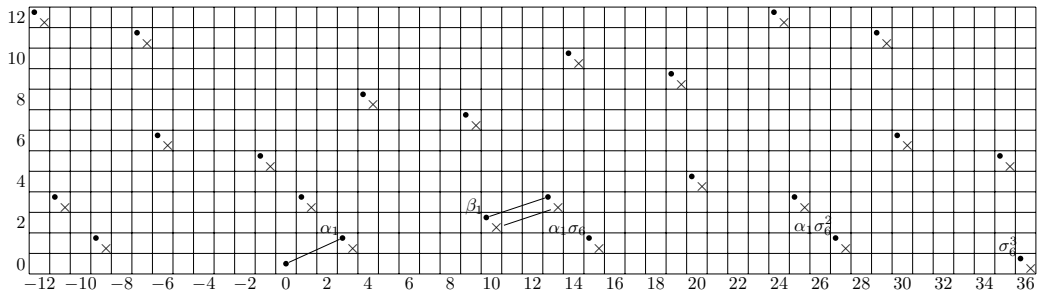}
\caption{The $E_6$-page of the homotopy fixed points}
\label{fig:HFPE6}
\end{figure}

The classes $\beta\U$ and $\U\dsix^2$, representing the classes
$\sqrt{-1}\beta_1$ and $\sqrt{-1}\dsix^3$ in the Tate spectral
sequence, are $d_5$-cycles. We depict the $E_6$-page of the homotopy
fixed point spectral sequence in Figure~\ref{fig:HFPE6}. In this
picture, both dots and crosses represent a copy of $\Z/3$. The
distinction reflects the action of the Atkin-Lehner involution
$t_5$. Classes denoted by a dot are fixed by $t_5$, while those
denoted by a cross are negated.  This involution commutes with the
differentials.

Sparseness of the spectral sequence guarantees that the next possible
differential is a $d_9$, and this is again governed by the Tate spectral
sequence.

\begin{cor}
There are $d_9$-differentials of the form
\[
d_9(\alpha_1\dsix\cdot \dsix^3)=\beta_1^5\text{ and
}d_9(\alpha_1\U\dsix^3)=\beta_1^5\U \dsix^{-1}.
\]
\end{cor}

\begin{figure}[h]
\includegraphics[angle=-90, width=\textwidth]{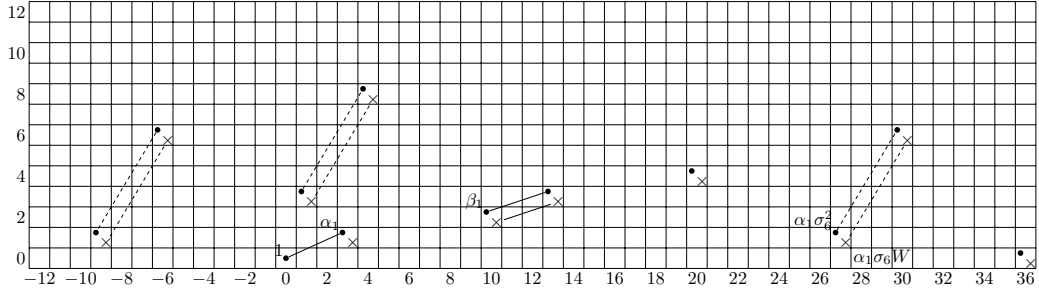}
\caption{The $E_\infty$-page of the homotopy fixed points}
\label{fig:HFPEinfty}
\end{figure}

These produce a horizontal vanishing line of $s$-intercept $10$, so there
are no further differentials possible. In particular, all elements in the
kernel of the map from the homotopy fixed point spectral sequence to the
Tate spectral sequence are permanent cycles. The $E_\infty$-page of the
spectral sequence (again ignoring classes on the zero-line in the
image of the transfer) is given in Figure~\ref{fig:HFPEinfty}.

There are a number of exotic multiplicative extensions. These are all
exactly analogous to those which arise with $\tmf$.

\begin{prop}
There are exotic $\alpha_1$-multiplications
\[
\alpha_1\cdot[\alpha_1\dsix^2]=\beta_1^3\text{ and
}\alpha_1\cdot[\alpha_1\dsix\U]=\beta_1^3\U \dsix^{-1}.
\]
\end{prop}

\begin{proof}
  The class $\dsix^2$, via the $d_5$-differential, represents a
  null-homotopy of $\alpha_1\beta_1^2$. By definition, we therefore
  conclude that
  $[\alpha_1\dsix^2]=\toda{\alpha_1,\alpha_1,\beta_1^2}.$ Standard
  shuffling results then show that
    \[
    \alpha_1\cdot\toda{\alpha_1,\alpha_1,\beta^2}=\toda{\alpha_1,\alpha_1,\alpha_1}\beta_1^2=\beta_1^3,
    \]
    giving the first result. The second is identical.
\end{proof}

In fact, these $\alpha_1$ extensions, together with the fact that
$\beta_1^2\alpha_1$ is zero, imply the $d_9$-differentials: $\dsix^4$ hits
$\beta_1^2\cdot\alpha_1\dsix^2$, so $\alpha_1\dsix^4$ hits
$\alpha_1\cdot\beta_1^2\alpha_1\dsix^2=\beta_1^5$. For degree reasons,
there are no other exotic multiplications, and there are no additive
extensions, so the computation is concluded.

The homotopy fixed point spectrum under the action of $t_2$ and $t_5$
is much more surprising.  Let $\taf^D_3$ be the connective cover of
$\TAF^D_3$.
\begin{thm}
As algebras, 
\[
\pi_\ast(\taf^D_3)^{h\Z/2\times\Z/2}\cong\pi_\ast\comp{\tmf}_3.
\]
\end{thm}

Our argument is purely computational. The analysis of the action of
$t_2$ and $t_5$ on the automorphic forms shows that the $E_2$-term for
the homotopy fixed point spectral sequence for ${\cal X}^\ast$ is
isomorphic to that of $\tmf$ \cite{bauer}. Upon identification of the
classes $\alpha_1$ and $\beta_1$, the arguments used are formal and
universal.  This leads us to the conclusion that the homotopy rings
are isomorphic, but for no obvious geometric reason.

\nocite{*}
\bibliography{shimc}

\begin{thebibliography}{25}
\expandafter\ifx\csname natexlab\endcsname\relax\def\natexlab#1{#1}\fi
\expandafter\ifx\csname url\endcsname\relax
  \def\url#1{\texttt{#1}}\fi
\expandafter\ifx\csname urlprefix\endcsname\relax\def\urlprefix{URL }\fi
\providecommand{\eprint}[2][]{\url{#2}}
\providecommand{\bibinfo}[2]{#2}
\ifx\xfnm\relax \def\xfnm[#1]{\unskip,\space#1}\fi
\bibitem[{Alsina and Bayer(2004)}]{alsinabayer}
\bibinfo{author}{Alsina, M.}, \bibinfo{author}{Bayer, P.},
  \bibinfo{year}{2004}.
\newblock \bibinfo{title}{Quaternion orders, quadratic forms, and {S}himura
  curves}. volume~\bibinfo{volume}{22} of \textit{\bibinfo{series}{CRM
  Monograph Series}}.
\newblock \bibinfo{publisher}{American Mathematical Society},
  \bibinfo{address}{Providence, RI}.
\bibitem[{Baba and Granath(2008)}]{babagranath}
\bibinfo{author}{Baba, S.}, \bibinfo{author}{Granath, H.},
  \bibinfo{year}{2008}.
\newblock \bibinfo{title}{Genus 2 curves with quaternionic multiplication}.
\newblock \bibinfo{journal}{Canad. J. Math.} \bibinfo{volume}{60},
  \bibinfo{pages}{734--757}.
\bibitem[{Baker(2000)}]{baker}
\bibinfo{author}{Baker, A.}, \bibinfo{year}{2000}.
\newblock \bibinfo{title}{{$I_n$}-local {J}ohnson-{W}ilson spectra and their
  {H}opf algebroids}.
\newblock \bibinfo{journal}{Doc. Math.} \bibinfo{volume}{5},
  \bibinfo{pages}{351--364 (electronic)}.
\bibitem[{Bauer(2008)}]{bauer}
\bibinfo{author}{Bauer, T.}, \bibinfo{year}{2008}.
\newblock \bibinfo{title}{Computation of the homotopy of the spectrum {$tmf$}},
  in: \bibinfo{booktitle}{Groups, homotopy and configuration spaces (Tokyo
  2005)}, \bibinfo{publisher}{Geom. Topol. Publ., Coventry}. pp.
  \bibinfo{pages}{11--40}.
\bibitem[{Behrens()}]{tmfconstr}
\bibinfo{author}{Behrens, M.}, .
\newblock \bibinfo{title}{Notes on the construction of {${\rm tmf}$}}.
\newblock \bibinfo{note}{Preprint.}
\bibitem[{Behrens and Lawson()}]{absurf}
\bibinfo{author}{Behrens, M.}, \bibinfo{author}{Lawson, T.}, .
\newblock \bibinfo{title}{Topological automorphic forms on ${U}(1,1)$}.
\newblock \bibinfo{note}{To appear in Math. Zeit.}
\bibitem[{Behrens and Lawson(2010)}]{taf}
\bibinfo{author}{Behrens, M.}, \bibinfo{author}{Lawson, T.},
  \bibinfo{year}{2010}.
\newblock \bibinfo{title}{Topological automorphic forms}.
\newblock \bibinfo{journal}{Mem. Amer. Math. Soc.} \bibinfo{volume}{204},
  \bibinfo{pages}{xxiii+136pp}.
\bibitem[{Boutot et~al.(1979)Boutot, Breen, G{\'e}rardin, Giraud, Labesse,
  Milne and Soul{\'e}}]{boutot}
\bibinfo{author}{Boutot, J.F.}, \bibinfo{author}{Breen, L.},
  \bibinfo{author}{G{\'e}rardin, P.}, \bibinfo{author}{Giraud, J.},
  \bibinfo{author}{Labesse, J.P.}, \bibinfo{author}{Milne, J.S.},
  \bibinfo{author}{Soul{\'e}, C.}, \bibinfo{year}{1979}.
\newblock \bibinfo{title}{Vari\'et\'es de {S}himura et fonctions {$L$}}.
  volume~\bibinfo{volume}{6} of \textit{\bibinfo{series}{Publications
  Math\'ematiques de l'Universit\'e Paris VII [Mathematical Publications of the
  University of Paris VII]}}.
\newblock \bibinfo{publisher}{Universit\'e de Paris VII U.E.R. de
  Math\'ematiques}, \bibinfo{address}{Paris}.
\bibitem[{Buzzard(1997)}]{buzzard}
\bibinfo{author}{Buzzard, K.}, \bibinfo{year}{1997}.
\newblock \bibinfo{title}{Integral models of certain {S}himura curves}.
\newblock \bibinfo{journal}{Duke Math. J.} \bibinfo{volume}{87},
  \bibinfo{pages}{591--612}.
\bibitem[{Deligne and Rapoport(1973)}]{delignerapoport}
\bibinfo{author}{Deligne, P.}, \bibinfo{author}{Rapoport, M.},
  \bibinfo{year}{1973}.
\newblock \bibinfo{title}{Les sch\'emas de modules de courbes elliptiques}, in:
  \bibinfo{booktitle}{Modular functions of one variable, {II} ({P}roc.
  {I}nternat. {S}ummer {S}chool, {U}niv. {A}ntwerp, {A}ntwerp, 1972)}.
  \bibinfo{publisher}{Springer}, \bibinfo{address}{Berlin}, pp.
  \bibinfo{pages}{143--316. Lecture Notes in Math., Vol. 349}.
\bibitem[{Devinatz et~al.(1988)Devinatz, Hopkins and Smith}]{nilpotence}
\bibinfo{author}{Devinatz, E.S.}, \bibinfo{author}{Hopkins, M.J.},
  \bibinfo{author}{Smith, J.H.}, \bibinfo{year}{1988}.
\newblock \bibinfo{title}{Nilpotence and stable homotopy theory. {I}}.
\newblock \bibinfo{journal}{Ann. of Math. (2)} \bibinfo{volume}{128},
  \bibinfo{pages}{207--241}.
\bibitem[{Eichler(1938)}]{eichler}
\bibinfo{author}{Eichler, M.}, \bibinfo{year}{1938}.
\newblock \bibinfo{title}{\"{U}ber die {I}dealklassenzahl hyperkomplexer
  {S}ysteme}.
\newblock \bibinfo{journal}{Math. Z.} \bibinfo{volume}{43},
  \bibinfo{pages}{481--494}.
\bibitem[{Elkies(1998)}]{elkies}
\bibinfo{author}{Elkies, N.D.}, \bibinfo{year}{1998}.
\newblock \bibinfo{title}{Shimura curve computations}, in:
  \bibinfo{booktitle}{Algorithmic number theory ({P}ortland, {OR}, 1998)}.
  \bibinfo{publisher}{Springer}, \bibinfo{address}{Berlin}. volume
  \bibinfo{volume}{1423} of \textit{\bibinfo{series}{Lecture Notes in Comput.
  Sci.}}, pp. \bibinfo{pages}{1--47}.
\bibitem[{Kudla et~al.(2006)Kudla, Rapoport and Yang}]{kry}
\bibinfo{author}{Kudla, S.S.}, \bibinfo{author}{Rapoport, M.},
  \bibinfo{author}{Yang, T.}, \bibinfo{year}{2006}.
\newblock \bibinfo{title}{Modular forms and special cycles on {S}himura
  curves}. volume \bibinfo{volume}{161} of \textit{\bibinfo{series}{Annals of
  Mathematics Studies}}.
\newblock \bibinfo{publisher}{Princeton University Press},
  \bibinfo{address}{Princeton, NJ}.
\bibitem[{Kurihara(1979)}]{kurihara}
\bibinfo{author}{Kurihara, A.}, \bibinfo{year}{1979}.
\newblock \bibinfo{title}{On some examples of equations defining {S}himura
  curves and the {M}umford uniformization}.
\newblock \bibinfo{journal}{J. Fac. Sci. Univ. Tokyo Sect. IA Math.}
  \bibinfo{volume}{25}, \bibinfo{pages}{277--300}.
\bibitem[{Milne(1979)}]{milne}
\bibinfo{author}{Milne, J.S.}, \bibinfo{year}{1979}.
\newblock \bibinfo{title}{Points on {S}himura varieties mod {$p$}}, in:
  \bibinfo{booktitle}{Automorphic forms, representations and {$L$}-functions
  ({P}roc. {S}ympos. {P}ure {M}ath., {O}regon {S}tate {U}niv., {C}orvallis,
  {O}re., 1977), {P}art 2}. \bibinfo{publisher}{Amer. Math. Soc.},
  \bibinfo{address}{Providence, R.I.}. Proc. Sympos. Pure Math., XXXIII, pp.
  \bibinfo{pages}{165--184}.
\bibitem[{Miyake(1989)}]{miyake}
\bibinfo{author}{Miyake, T.}, \bibinfo{year}{1989}.
\newblock \bibinfo{title}{Modular forms}.
\newblock \bibinfo{publisher}{Springer-Verlag}, \bibinfo{address}{Berlin}.
\newblock \bibinfo{note}{Translated from the Japanese by Yoshitaka Maeda}.
\bibitem[{Morita(1981)}]{morita}
\bibinfo{author}{Morita, Y.}, \bibinfo{year}{1981}.
\newblock \bibinfo{title}{Reduction modulo {${P}$} of {S}himura curves}.
\newblock \bibinfo{journal}{Hokkaido Math. J.} \bibinfo{volume}{10},
  \bibinfo{pages}{209--238}.
\bibitem[{Mumford(1970)}]{mumford}
\bibinfo{author}{Mumford, D.}, \bibinfo{year}{1970}.
\newblock \bibinfo{title}{Abelian varieties}.
\newblock Tata Institute of Fundamental Research Studies in Mathematics, No. 5,
  \bibinfo{publisher}{Published for the Tata Institute of Fundamental Research,
  Bombay}.
\bibitem[{Oort(1971)}]{oort}
\bibinfo{author}{Oort, F.}, \bibinfo{year}{1971}.
\newblock \bibinfo{title}{Finite group schemes, local moduli for abelian
  varieties, and lifting problems}.
\newblock \bibinfo{journal}{Compositio Math.} \bibinfo{volume}{23},
  \bibinfo{pages}{265--296}.
\bibitem[{Ravenel(1986)}]{Greenbook}
\bibinfo{author}{Ravenel, D.C.}, \bibinfo{year}{1986}.
\newblock \bibinfo{title}{Complex cobordism and stable homotopy groups of
  spheres}. volume \bibinfo{volume}{121} of \textit{\bibinfo{series}{Pure and
  Applied Mathematics}}.
\newblock \bibinfo{publisher}{Academic Press Inc.}, \bibinfo{address}{Orlando,
  FL}.
\bibitem[{Serre(1973)}]{serrearith}
\bibinfo{author}{Serre, J.P.}, \bibinfo{year}{1973}.
\newblock \bibinfo{title}{A course in arithmetic}.
\newblock \bibinfo{publisher}{Springer-Verlag}, \bibinfo{address}{New York}.
\newblock \bibinfo{note}{Translated from the French, Graduate Texts in
  Mathematics, No. 7}.
\bibitem[{Shimura(1967)}]{shimura}
\bibinfo{author}{Shimura, G.}, \bibinfo{year}{1967}.
\newblock \bibinfo{title}{Construction of class fields and zeta functions of
  algebraic curves}.
\newblock \bibinfo{journal}{Ann. of Math. (2)} \bibinfo{volume}{85},
  \bibinfo{pages}{58--159}.
\bibitem[{Shimura(1975)}]{shimurareal}
\bibinfo{author}{Shimura, G.}, \bibinfo{year}{1975}.
\newblock \bibinfo{title}{On the real points of an arithmetic quotient of a
  bounded symmetric domain}.
\newblock \bibinfo{journal}{Math. Ann.} \bibinfo{volume}{215},
  \bibinfo{pages}{135--164}.
\bibitem[{Strickland(1999)}]{strickland}
\bibinfo{author}{Strickland, N.P.}, \bibinfo{year}{1999}.
\newblock \bibinfo{title}{Products on {${\rm MU}$}-modules}.
\newblock \bibinfo{journal}{Trans. Amer. Math. Soc.} \bibinfo{volume}{351},
  \bibinfo{pages}{2569--2606}.

\end{thebibliography}

\end{document}